[1994/12/01]% LaTeX date must December 1994 or later
\documentclass[11pt,reqno]{article}

\title{Recurrence for branching Markov chains}
\author{Sebastian M\"uller\thanks{partially supported by FWF (Austrian Science Fund) project P19115-N18
and DFG (German Research Foundation) project MU 2868/1-1} \\
Institut f\"ur mathematische Strukturtheorie\\
Technische Universit\"at Graz\\
Steyrergasse 30/III  \\
8010 Graz\\
Austria\\
{\tt mueller@tugraz.at}\\}
\usepackage{amsmath,amsthm,amsfonts,bm}
\usepackage[latin1]{inputenc}
\usepackage{t1enc}
\usepackage{pst-all}
\usepackage{amsmath}        % use AMS packages
\usepackage{amssymb}
\usepackage{amsthm}

\newtheorem{thm}{Theorem}[section]
\newtheorem{cor}[thm]{Corollary}
\newtheorem{lem}[thm]{Lemma}
\newtheorem{prop}[thm]{Proposition}
\newtheorem{conj}[thm]{Conjecture}

\theoremstyle{definition}
\newtheorem{defn}[thm]{Definition}
\theoremstyle{remark}
\newtheorem{rem}[thm]{Remark}

\newtheorem{ex}[thm]{Example}
\numberwithin{equation}{section}

\def\orig{o}
\newcommand{\eps}{\varepsilon}

\newcommand{\R}{\mathbb{R}}
\newcommand{\Z}{\mathbb{Z}}
\newcommand{\N}{\mathbb{N}}
\renewcommand{\P}{\mathbb{P}}
\newcommand{\E}{\mathbb{E}}

\newcommand{\kommentar}[1]{}

\def\orig{o}

\def\0b{\bf{0}}
\def\R{\mathbb{R}}
\def\T{\mathbb{T}}
\def\P{\mathbb{P}}
\def\E{\mathbb{E}}

\def\Z{\mathbb{Z}}

\def\N{\mathbb{N}}

\def\eps{\varepsilon}
\def\vr{\varrho}

\begin{document}
\maketitle {\abstract The question of recurrence and transience of
branching Markov chains is more subtle than for ordinary Markov
chains; they can be  classified in transience, weak recurrence,
and strong recurrence. We review criteria for transience and weak
recurrence and give several new conditions for weak recurrence and
strong recurrence.  These conditions make a unified treatment of
known and new examples possible and provide enough information to
distinguish between weak and strong recurrence. This represents a
step towards a general classification of branching Markov chains.
In particular, we show that in \emph{homogeneous} cases weak
recurrence and strong recurrence coincide. Furthermore, we discuss
the generalization of positive and null recurrence to branching
Markov chains and show that branching random walks on $\Z$ are
either transient or positive recurrent.
\renewcommand{\sectionmark}[1]{}
\newpage

\section{Introduction}

A branching Markov chain (BMC) $(X,P,\mu)$  is  a system of
particles in discrete time on a discrete state space $X.$ The
process starts with one particle in some starting position $x\in
X.$ At each time particles split up in offspring particles
independently according to some probability distributions
$(\mu(x))_{x\in X}.$ The new particles then move independently
according to some irreducible Markov chain $(X,P)$ with transition
probabilities $P.$ Processes of this type are studied in various
articles with different notations and variations of the model.

%One aim of this paper is to review the existing results and to
%propose a unifying description. Furthermore, we present new
%condition for strong rec

Ordinary  Markov chains are either transient or recurrent, i.e.,
the starting position is either visited a finite or an infinite
number of times.  This $0-1$-law does not hold for branching
Markov chains in general, compare with \cite{benjamini:94},
\cite{gantert:04}, and \cite{mueller06}. Let $\alpha(x)$ be the
probability that starting the BMC in $x$ the state $x$ is visited
an infinite number of times by some particles. We can classify the
BMC in transient ($\alpha(x)=0$ for all $x$), weakly recurrent
($0<\alpha(x)<1$ for all $x$), and strongly recurrent
($\alpha(x)=1$ for all $x$). In cases where we do not want to
distinguish between weak and strong recurrence we just say that
the process is recurrent. Effects of this type occur also in a
varied model of BMC in which the branching is governed by a fixed
tree, compare with \cite{benjamini:94} for more details on Markov
chains indexed by trees.

Let $m(x)$ be the mean of the offspring distribution $\mu(x).$ If
$m(x)=m$ for all $x$ there is a nice and general criterion for
transience and recurrence: the BMC $(X,P,\mu)$ is transient if and
only if $m\leq 1/\rho(P),$ where $\rho(P)$ is the spectral radius
of the underlying Markov chain $(X, P).$ Observe that the type of
the process depends only on two characteristics, $m$ and
$\rho(P),$ of the process. We give several description of the
spectral radius $\rho(P).$  These description are useful to give
several different proofs for the above criterion and offer various
possibilities to decide whether a process is transient or
recurrent.

Another purpose of this paper is to review and to continue the
work of \cite{benjamini:94}, \cite{hollander:95},
\cite{gantert:04}, \cite{hueter:00}, \cite{machado:00},
\cite{machado:03}, \cite{menshikov:97}, \cite{mueller07},
\cite{mueller06}, \cite{pemantle:01}, \cite{schinazi:93},
\cite{stacey:03}, and \cite{volkov:01}. We present a unifying
description and give new  conditions for weak and strong
recurrence. Our results suggest that the type of the process only
depends on the mean number of offspring and some variation of the
spectral radius of the underlying Markov chain, too.

A useful observation concerning the different types is that
transience / recurrence is sensitive to local changes of the
underlying process but criteria for strong recurrence have to
contain  global information of the process. This fact is exploited
by the idea of \emph{seeds,} finite regions of the process that
may produce an infinite number of particles independently of the
rest of the process. Eventually, infinitely many particles leave
this finite regions and return to the starting position. Hence,
the BMC is recurrent. While the existence of one such seed already
implies recurrence, we need \emph{sufficiently many good} seeds to
guarantee strong recurrence, compare with Section
\ref{sec:strong}. In \emph{homogeneous} cases, where local
properties become global, recurrence and strong recurrence
coincide, compare with Subsection \ref{subsec:homo}.

In Section \ref{sec:pre} we recall some facts about Markov chains,
Green functions, and the corresponding spectral radius that are
crucial for our further development. In particular, we give
several description of the spectral radius of a Markov chain, e.g.
in terms of power series, superharmonic functions,
Perron-Frobenius eigenvalues, and rate function of large
deviations.  In Subsection \ref{subsec:bmc} we define the model of
BMC formally, recall the classification results  for transience
and recurrence of \cite{gantert:04} and \cite{mueller06}, and give
concrete examples in  Subsection \ref{subsec:ex}.

In Section \ref{sec:strong} we first give general conditions for
strong recurrence, see Subsection \ref{sec:gen_crit}, and  show
that in \emph{homogeneous} cases recurrence and strong recurrence
coincide, see Subsection \ref{subsec:homo} where we prove known
(Theorems \ref{thm:2} and \ref{thm:5}) and new results (Theorems
\ref{thm:2_finitecone} and \ref{thm:7} as well as Subsection
\ref{sec:uniform}) with a unified method. We then discuss, see
Subsections \ref{subsect:finer}, several conditions for strong
recurrence that, while not leading to a complete classification in
general, work well in concrete situations and offer a toolbox to
handle further concrete problems. In particular, there is a second
critical value that divides the weakly and strongly recurrent
phase. We want to point out that at the second critical value the
process may be either strongly or weakly recurrent, compare with
Theorem \ref{thm:rec_quasi}. In order to develop a complete
classification we study BMC on graphs in Subsection
\ref{subsec:rwgraphs} and give conjectures in Subsection
\ref{sec:outlook}.

In Section \ref{sec:posrec} we generalize the concept of positive
and null recurrence to BMC and give conditions for positive
recurrence in terms of a functional equation, Theorem
\ref{thm:bmc:2Foster}, that fits in the context of Theorems
\ref{thm:1} and \ref{thm:strecur} and can be seen as the natural
generalization of the criterion of Foster for positive recurrence
of Markov chains. Eventually, we conclude in showing that
homogeneous branching Markov chains on $\Z$ with bounded jumps are
either transient or positive recurrent, compare with Theorem
\ref{thm:posrecZ}.

\section{Preliminaries}\label{sec:pre}

\subsection{Markov chains}
A  Markov chain  $(X,P)$ is defined by a countable  state space
$X$ and  transition probabilities $P=\left( p(x,y) \right)_{x,y\in
X}.$ The elements $p(x,y)$ of $P$ define the probability of moving
from $x$ to $y$ in one step. As usual, let $X_n$ denote the
position of the Markov chain at time $n.$ The transition operator
$P$ can be interpreted as a (countable) stochastic matrix, so
that, on the one hand, $p^{(n)}(x,y)$ is the $(x,y)$-entry of the
matrix power $P^n$ and, on the other hand, we have that
$p^{(n)}(x,y)=\P_x(X_n=y)$ is the probability to get from $x$ to
$y$ in $n$ steps. We set $P^0=I,$ the identity matrix over $X.$
Throughout this paper we assume that the state space is infinite
and that the transition operator $P$ is irreducible, i.e., for
every pair $x,y\in X$ there exists some $k\in\N$ such that
$p^{(k)}(x,y)>0.$

Markov chains are related to  random walks on graphs.  So let us
recall some basic standard notations of graphs. A (directed) graph
$G=(V,E)$ consists of a finite or countable set of vertices $V$
and an adjacency relation $\sim$ that defines the set of edges
$E\subset V\times V.$ A   path from a vertex $x$ to some vertex
$y$ is a sequence $x=x_0,x_1,\ldots,x_n=y$ with $x_i\sim x_{i+1}$
for all $0\leq i< n.$ The number $n$ is the length of the path. A
graph is (strongly) connected if every ordered pair of vertices is
joined by a path. The usual graph distance $d(x,y)$ is the minimum
among the length of all paths from $x$ to $y.$ A vertex $y$ is
called  a neighbor of $x$ if $x\sim y.$ The degree $deg(x)$ of a
vertex $x$ is the numbers of its neighbors. A graph is called
locally finite if $deg(x)$ is finite for all vertices $x$. We say
a graph $G$ has  bounded geometry if $deg(\cdot)$ is bounded and
is $M$-regular if all vertices have degree $M.$

Every Markov chain $(X,P)$  defines a  graph $G=(V,E),$ with a set
of vertices $V=X$ and a set of edges $E:=\{(x,y): p(x,y)>0 ~x,y\in
X\}.$ It is clear that a Markov chain is irreducible if and only
if its corresponding graph is connected. If the transition
probabilities of the Markov chain are in some kind adapted to the
structure of $G,$ we shall speak of a random walk on $G$ with
transition probabilities $P.$  We shall call a Markov chain on a
graph with symmetric adjacency relation a simple random walk
(SRW) if the {\it walker} chooses every neighbor  with the same
probability, i.e., $p(x,y)=1/deg(x)$ for $x\sim y$ and $0$
otherwise.
%The graph distance between $x$ and $y$ is the minimum
%number of steps needed  to get from $x$ to $y:$
%$$ d(x,y)=\inf_n\{p^{(n)}(x,y)>0\}.$$

We recall the Green function and the spectral radius of an
irreducible Markov chain. These two characteristics will be
crucial for our further development, compare with \S 1 in
\cite{woess} for proofs and more.

\begin{defn}\label{def:green_fct}
The Green function of $(X,P)$ is the power series
$$G(x,y|z):=\sum_{n=0}^\infty p^{(n)} (x,y) z^n,~ x,y\in X,~
z\in\mathbb{C}.$$ We write $G(x,x)$ for $G(x,x|1).$
\end{defn}

Observe, that, due to the irreducibility, $G(x,y|z)$ either
converges for all $x,y\in X$ or diverges for all $x,y\in X.$
Therefore, we can define  $R$ as the finite convergence radius of
the series $G(x,y|z)$ and call $1/R$ the spectral radius of the
Markov chain.

\begin{defn}\label{def:spectr_rad}
 The spectral  radius of $(X,P)$ is defined as
\begin{equation}\label{def:spectr}
\rho(P):=\limsup_{n\rightarrow\infty}
\left(p^{(n)}(x,y)\right)^{1/n} \in (0,1].
\end{equation}
We denote $\rho(G)$ the spectral radius of the SRW on the graph
$G.$
\end{defn}

\begin{lem}\label{lem:spec_sup} We have
$$p^{(n)}(x,x)\leq \rho(P)^n,~and ~
\lim_{n\rightarrow\infty}\left(p^{(nd)}(x,x)\right)^{1/nd}=\rho(P),$$
where  $d:=d(P):=gcd\{n:~ p^{(n)}(x,x)>0~\forall x\}$ is the {\it
period} of $P.$
\end{lem}

If $X$ and $P$  are finite, the spectral radius of $\rho(P)$
becomes the largest, in absolute value, eigenvalue of the matrix
$P$ and equals $1.$ The spectral radius of a Markov chain with
infinite state space can be approximated with  spectral radii of
finite sub-stochastic matrices. To this end we  consider (general)
finite nonnegative matrices $Q.$  A matrix
$Q=(Q_{i,j})_{i,j\in\R^{N\times N}}$ with nonnegative entries is
called  irreducible if for any pair of indices $i,j$ we have
$Q^m(i,j)>0$ for some $m\in \N.$ The well-known Perron-Frobenius
 Theorem states (e.g. Theorem 3.1.1 in \cite{dembo92}), among other things, that $Q$
possesses a largest, in absolute value, real eigenvalue $\rho(Q).$
Furthermore, we have (e.g. with part \emph{(e)} of Theorem 3.1.1
in \cite{dembo92})
$$\rho(Q)=\limsup_{n\to\infty}\left(Q^n(i,i)\right)^{1/n}\quad
\forall 1\leq i\leq N.$$

\begin{rem}
If $Q$ is a symmetric  matrix with $\rho(Q)<\infty$ then $Q$ acts
on $l^2(X)$ as a bounded linear operator with norm
$\|Q\|=\rho(Q).$ The same holds true for reversible Markov chains
$(X,P)$ with $P$ acting on some appropriate Hilbert space.
\end{rem}

Now, let us consider  an infinite irreducible Markov chain $(X,P).$ A
subset $Y\subset X$ is called irreducible if the sub-stochastic
operator $$P_Y=(p_Y(x,y))_{x,y\in Y}$$ defined by
$p_Y(x,y):=p(x,y)$ for all $x,y\in Y$ is irreducible.

It is straightforward to show the next characterization.

\begin{lem}\label{specapprox}
Let $(X, P)$ be an irreducible Markov chain. Then, \begin{equation}\label{eq:specapprox}
\rho(P)=\sup_Y \rho(P_Y),
\end{equation}
where the supremum is over finite and irreducible subsets
$Y\subset X.$ Furthermore, $\rho(P_F)< \rho(P_G)$ if $F\subset
G.$
\end{lem}

For finite irreducible matrices $Q$ the spectral radius as defined
in (\ref{def:spectr}) equals the largest eigenvalue of $Q.$ This
does not longer hold true for general infinite irreducible
transition kernels $P,$ compare with \cite{vere-jones:67}, since
the existence of a Perron-Frobenius eigenvalue can not be
guaranteed. Nevertheless, the transition operator $P$ acts on
functions $f:~X\rightarrow\R$ by
\begin{equation}\label{eq:Pf}
Pf(x):= \sum_y p(x,y) f(y),
\end{equation}
where we assume that $P|f|<\infty.$ It turns out that the
spectral radius $\rho(P)$  can be characterized in terms of $t$-superharmonic  functions.

\begin{defn} Fix $t>0.$
A  $t$-superharmonic function  is a function
$f:~X\rightarrow\R$ satisfying
$$Pf(x)\leq t\cdot f(x) \quad \forall x\in X.$$
\end{defn}

We obtain the following (well-)known characterization of the
spectral radius in terms of $t$-superharmonic functions, e.g.
compare with  \S 7 in \cite{woess}.
\begin{lem}\label{lem:1}
$$\rho(P)=\min\{t>0:~\exists\, f(\cdot)>0\mbox{ such that }Pf\leq t f\}$$
\end{lem}

We can express the spectral radius in terms of the rate function $I(\cdot)$ of a large deviation principle (LDP).
Let us assume that a LDP holds for the distance, i.e.,

\begin{eqnarray*}
  \liminf_{n\to \infty} \frac1n \log
  \P_\orig\left(\frac{d(X_n,\orig)}n\in  O\right) &\geq& -\inf_{a\in O} I(a), \\
 \limsup_{n\to \infty} \frac1n \log
 \P_\orig\left(\frac{d(X_n,\orig)}n\in  C\right) &\leq& -\inf_{a\in C} I(a),
\end{eqnarray*}
for all open sets $O\subset \R$ and all closed sets $C\subset \R.$
Furthermore, we assume ellipticity, i.e.,  that
 there is  some constant $c>0$ such that $p(x,y)\geq c\cdot
\mathbf{1}\{p(x,y)>0\}$ for all $x,y\in X.$

We do not claim to be the first who make the following observation
that is quite intuitive since LDP's for Markov chains are closely
linked with Perron-Frobenius theory for irreducible matrices.

\begin{lem}\label{lem:spectral_ldp}
Let $(X,P)$ be an irreducible Markov chain with $p(x,y)\geq c\cdot
\mathbf{1}\{p(x,y)>0\}$ for some constant $c>0.$ Assume
furthermore that a LDP holds  for the distance with rate function
$I(\cdot),$ then
\begin{equation}
-\log \rho(P)= I(0).
\end{equation}
\end{lem}
\begin{proof}
We have
\begin{eqnarray*}
  - I(0) &=& \lim_{\eps\rightarrow 0} \lim_{n\rightarrow \infty}
  \frac1n \log \P_\orig \left(\frac{d(X_n,\orig)}n\leq \eps\right) \\
  &\geq &\limsup_{n\to\infty} \frac1n \log \P_\orig\left(X_n=\orig\right)=
  \log \rho(P).
\end{eqnarray*} For the converse we use our assumption: $p(x,y)\geq c\cdot
\mathbf{1}\{p(x,y)>0\}$ for some $c>0.$ With $0<\eps <1$ we obtain

$$p^{(n)}(\orig,\orig)\geq \P_\orig\left(d(X_{\lceil n(1-\eps)\rceil},\orig)
\leq  \lfloor n\eps\rfloor\right) c^{\lfloor n\eps\rfloor}.$$
Therefore,
\begin{eqnarray*}
  \frac 1n \log \P_\orig\left( \frac{d(X_{\lceil n(1-\eps)\rceil},\orig)}
  {\lceil n(1-\eps)\rceil} \leq
  {\frac{\lfloor n\eps\rfloor}{\lceil n(1-\eps)\rceil}}\right)
   &=& \frac 1n \log \P_\orig\left(d(X_{\lceil n(1-\eps)\rceil},
   \orig)\leq  \lfloor n\eps\rfloor\right) \\
   &\leq & \frac1n \log \left(p^{(n)}(\orig,\orig) c^{-\lfloor
   n\eps\rfloor})\right) \\
   &=&
   \log \left(p^{(n)}(\orig,\orig)\right)^{1/n}-
   \frac{\lfloor n\eps\rfloor}{n} \log c.
\end{eqnarray*}
Hence, for all $0< \eps <1$ $$ -I(0)\leq \lim_{n\rightarrow\infty}
\frac1n \log \P_\orig\left(
\frac{d(X_{n},\orig)}{n}\leq{\frac{\lfloor n\eps\rfloor}{\lceil
n(1-\eps)\rceil}}\right)\leq \log\rho - \eps\log c.$$
Now, letting $\eps\to0$ finishes the proof.
\end{proof}

\subsection{Branching Markov chains}\label{subsec:bmc}

A branching Markov chain (BMC) $(X,P,\mu)$ consists of two
dynamics on the state space $X$: branching, $(\mu(x))_{x\in X},$
and movement, $P.$ Here
$$\mu(x)=\left(\mu_k(x)\right)_{k\geq 1}$$ is a sequence of
nonnegative numbers
satisfying
$$\sum_{k=1}^\infty \mu_k(x)=1 \mbox{ and }
m(x):=\sum_{k=1}^\infty k \mu_k(x)<\infty,$$ and $\mu_k(x)$ is the
probability that a particle in $x$ splits up in $k$ particles. The
movement of the particles is governed through an irreducible and
infinite Markov chain $(X,P).$ Note that we always  assume that
each particle has at least one offspring, i.e., $\mu_0=0,$ and
hence the process survives forever. Similar results can be
obtained by conditioning on the survival of the branching process.

The BMC  is defined as follows. At time $0$ we start with one
particle in an arbitrary starting position $x\in X.$  At time $1$
this particle splits up in $k$ offspring particles  with
probability $\mu_k(x).$ Still at time $n=1,$ these $k$ offspring
particles then  move independently according to the Markov chain
$(X,P).$ The process is defined inductively: At each time each
particle in position $x$ splits up according to $\mu(x)$ and the
offspring particles move according to $(X,P).$ At any time, all
particles move and branch independently of the other particles and
the previous history of the process. If the underlying Markov
chain is a random walk on a graph $G$ and the branching
distributions are adapted to $G$ we shall also speak of a
branching random walk (BRW). In the case where the underlying
random walk is a simple random walk on a graph $G$ we denote the
process $(G,\mu).$

We introduce the following notations in order to describe the
process. Let $\eta(n)$ be the total number of particles at time
$n$ and let $x_i(n)$ denote the position of the $i$th particle at
time $n.$ Denote $\P_x(\cdot)=\P(\cdot|x_1(0)=x)$ the probability
measure for  a BMC started with one particle in $x.$ A BMC can be
viewed as a Markov chain on the \emph{big} state space
$$\mathcal{X}:= X \cup X^2 \cup X^3 \cup \ldots~ .$$  Clearly, this
Markov chain on $\mathcal{X}$ is transient if $m(x)>1$ for some
$x\in X.$ A priori it is not clear in which sense transience and
recurrence of Markov chains can be generalized to BMC. One
possibility is to say a BMC is recurrent if with probability $1$
at least one particle returns to the starting position. This
approach was followed for example in \cite{menshikov:97}. We
choose a different one, e.g. compare with \cite{benjamini:94},
\cite{gantert:04}, and \cite{mueller06}, since this approach
offers a finer partition in recurrence and transience, and gives
interesting generalizations of equivalent conditions for
transience of (nonbranching) Markov chains, compare with Theorem
\ref{thm:1}.

\begin{defn}\label{def:bmc_rec} Let
\begin{equation}\label{eq:1} \alpha(x):=\P_x\left(
\sum_{n=1}^\infty
\sum_{i=1}^{\eta(n)}\mathbf{1}\{x_i(n)=x\}=\infty\right).
\end{equation}
be the probability that $x$ is visited an infinite number of
times, given the process starts in $x$ with one particle.
 A BMC is called transient if $\alpha(x)=0$ for some ($\Leftrightarrow$
all)  $x\in X,$ weakly recurrent if $0<\alpha(x)<1$ for some
($\Leftrightarrow$ all)  $x\in X,$ and
  strongly recurrent if $\alpha(x)=1$ for some ($\Leftrightarrow$
all) $x\in X.$   We write $\alpha=0$ if $\alpha(x)=0$ for all
$x\in X$ and $\alpha>0$ and $\alpha\equiv 1,$ respectively. We
call a BMC recurrent if it is not transient.
\end{defn}

Definition \ref{def:bmc_rec} is justified through the next result.

\begin{lem}\label{lem:bmc_rec}
We have that $\alpha(x)=0,~ \alpha(x)>0,$ and $\alpha(x)=1$ either
hold for all or none $x\in X.$
\end{lem}

\begin{proof}
Let $x,y\in X.$ Due to irreducibility and the independence of
branching and movement we have that $$\P_x\left(\sum_{n=1}^\infty
\sum_{i=1}^{\eta(n)}\mathbf{1}\{x_i(n)=y\}=\infty \right)
\begin{array}{c}
  = 0 \\
  >0  \\
  =1 \\
\end{array}$$ is equivalent to $$ \P_y\left(\sum_{n=1}^\infty
\sum_{i=1}^{\eta(n)}\mathbf{1}\{x_i(n)=y\}=\infty \right)
\begin{array}{c}
  = 0 \\
  >0  \\
  =1 \\
\end{array}. $$ Hence, it suffices to show that
\begin{equation}\label{eq:thm:1}\P_{x}\left( \sum_{n=1}^\infty
\sum_{i=1}^{\eta(n)}\mathbf{1}\{x_i(n)=x\}=\infty\mbox{ and }
 \sum_{n=1}^\infty
\sum_{i=1}^{\eta(n)}\mathbf{1}\{x_i(n)=y\}<\infty  \right)=0.
\end{equation} We follow the line of the
proof of Lemma 3.3 in \cite{benjamini:94}. Since $(X,P)$ is
irreducible we have $p^{(l)}(x,y)=\delta>0$ for some $l\in \N.$
Let $N, M\in\N.$ The probability that there are times
$M<n_1,\ldots,n_N$ with $n_{j-1}+l<n_j$ for all $1\leq j\leq N$
such that $x_i(n_j)=x$ for some $1\leq i\leq \eta(n_j)$ for all
$j$ but $x_i(n)\neq y$ for all $n>M$ and all $1\leq i \leq
\eta(n)$ is at most $(1-\delta)^N.$  Letting $N\rightarrow\infty,$
this yields
$$\P_{x}\left( \sum_{n=1}^\infty
\sum_{i=1}^{\eta(n)}\mathbf{1}\{x_i(n)=x\}=\infty\mbox{ and }
\sum_{n=M}^\infty
\sum_{i=1}^{\eta(n)}\mathbf{1}\{x_i(n)=y\}=0\right)=0.$$ Let $A_M$
be the event in the last formula. Notice that
$$\bigcup_{M\geq 1} A_M =\left\{ \sum_{n=1}^\infty
\sum_{i=1}^{\eta(n)}\mathbf{1}\{x_i(n)=x\}=\infty\mbox{ and }
\sum_{n=1}^\infty
\sum_{i=1}^{\eta(n)}\mathbf{1}\{x_i(n)=y\}<\infty\right\}.$$ This
proves the claim.
\end{proof}

In analogy to \cite{menshikov:97}, see also the superposition
principle in \cite{hueter:00}, we introduce the following modified
version of BMC. We fix some position $\orig\in X,$ which we denote
the origin of $X.$  The new process is like the original BMC at
time $n=1,$ but is different for $n>1.$ After the first time step
we conceive the origin as {\it freezing}: if a particle reaches
the origin, it stays there forever and stops splitting up. We
denote this new process with {\bf BMC*}. The process BMC* is
analogous to the original process BMC except that
$p(\orig,\orig)=1,~ p(\orig,x)=0~\forall x\neq \orig$ and
$\mu_1(\orig)=1$ from the second time step on. Let $\eta(n,x)$ be
the number of particles at position $x$ at time $n.$ We define the
random variable $\nu(\orig)$ as
$$\nu(\orig):=\lim_{n\rightarrow\infty} \eta(n,\orig)\in \{0,1,\ldots\}\cup\{\infty\}.$$
 We write $\E_x \nu(\orig)$ for the
expectation of $\nu(\orig)$ given that $x_1(0)=x.$  Note that our
notation of $\eta(n,\orig)$ and $\nu(\orig)$ is different from the
one in \cite{gantert:04}. Since the choice of the origin may
affect the behavior of the BMC*, we keep track of the dependence
of the variables $\eta$ and $\nu$ on the choice of the origin and
write $\eta(n,\orig)$ and $\nu(\orig).$ Furthermore, our
definition of the process BMC* differs from the one given in
\cite{menshikov:97}: in our definition, the origin is not
absorbing at time $n=1.$

The Green function $G(x,x|m)$ at $z=m$ gives the expected number
of particles that visits $x$ of a BMC with constant mean offspring
$m=m(x)$ started in $x.$ Due to this interpretation,
$G(x,y|m)<\infty$ implies transience of a BMC with constant mean
offspring $m.$ While in some cases, e.g. random walk on
homogeneous trees \cite{hueter:00}, the converse is  true, it does
not hold in general, compare with Remark \ref{rem:green}. It turns
out that another generating function is decisive for transience.
Let
$$ T_y:=\min_{n\geq 1}\{X_n=y\}$$ be the time of the first return to $y$
and $$U(x,y|z):=\sum_{n=1}^\infty \P_x(T_y=n) z^n$$ its
corresponding generating function. Due to the definition of the
process BMC* we have the useful identity
$$ \E_x \nu(y) = U(x,y|m)$$ for a BMC with constant mean offspring
$m.$ Furthermore, the equality offers a probabilistic
interpretation for the generating function $U(x,y|z)$ with $z\geq
1.$ The generating functions $G$ and $U$ are naturally connected
through
$$G(x,x|z)=\frac1{1-U(x,x|z)}$$ and hence one can show that

\begin{equation}\label{eq:rho_U}
\rho(P)=\max\{z>0:~ U(x,x|z)\leq 1\}.
\end{equation}

The next Theorem, due to \cite{mueller06}, gives several
sufficient and necessary conditions for transience. Notice that
the expected number of particles visiting $x$ the first time in
their ancestry line, i.e., $\E_x \nu(x),$ takes the role of the
Green function in the theory of Markov chains and that the
criterion of transience in terms of the existence of nonconstant
superharmonic functions becomes $(iii).$

\begin{thm}\label{thm:1}
A BMC $(X,P,\mu)$ with $m(y)>1$ for some $y$ is transient if and
only if the three equivalent conditions hold:
\begin{enumerate}
     \item[(i)] $\E_{\orig} \nu(\orig)\leq 1$ for some
     ($\Leftrightarrow$ all) $\orig\in X.$
     \item[(ii)] $\E_x \nu(\orig)<\infty$ for all $x,\orig\in X.$
     \item[(iii)] There
exists a strictly positive function $f(\cdot)$ such that
\begin{equation*}%\label{eq:1.1}
P f(x)\leq \frac{f(x)}{m(x)}\quad \forall x\in X.
\end{equation*}

\end{enumerate}
\end{thm}

In particular if the mean offspring is constant, i.e., $m(x)=m
~\forall x\in X,$ we have, due to  Lemma \ref{lem:1}, the
following  result of \cite{gantert:04}. Observe that in this case
we can speak about a \emph{critical behaviour.}

\begin{thm}\label{thm:1a}
The BMC  $(X,P,\mu)$ with constant mean offspring $m>1$ is
transient if  and only if $m\leq 1/\rho(P).$
\end{thm}
Theorem \ref{thm:1a} follows directly from Lemma \ref{lem:1} and
part \emph{(iii)} of  Theorem \ref{thm:1}. Another way, compare
with \cite{woess2}, to see Theorem \ref{thm:1a} is combining
Theorem \ref{thm:1} \emph{(i)} and the fact that
$\rho(P)=\max\{z>0:~ U(x,x|z)\leq 1\}$ (see Equation
(\ref{eq:rho_U})) and conclude with $\E_x \nu(y) = U(x,y|m).$ We
give a direct proof, without using the abstract arguments of Lemma
\ref{lem:1} or Equation (\ref{eq:rho_U}), since the arguments used
are helpful to understand the reasonings in the remainder of the
paper.

\begin{proof}
The case $m<1/\rho(P)$ is clear, since $G(x,x|m)<\infty$ implies
transience. To show that  $m >1/\rho(P)$ implies recurrence we
compare the original BMC with an embedded process and prove that
this process with fewer particles is recurrent. We start the BMC
in $\orig\in X$. We know from the hypothesis and the definition of
$\rho(P)$ that there exists a $k=k(\orig)$ such that
$$p^{(k)}(\orig,\orig)> m^{-k}.$$ We construct the embedded process
$(\xi_i)_{i\geq 0}$ by observing the BMC only  at times $k, 2k,
3k, \ldots $ and by neglecting all the particles not being in
position $\orig$ at these times. Let $\xi_i$ be the number of
particles of the new process in $\orig$ at time $ik.$ The process
$(\xi_i)_{i\geq 0}$ is a Galton-Watson process with mean
$p^{(k)}(\orig,\orig)\cdot m^{k} > 1,$ thus survives with positive
probability. Eventually, the origin is hit infinitely often with
positive probability.

In order to prove  transience at  the critical value
$m=1/\rho(P),$ we use a continuity argument to show that the
subset $\{m:~ (X,P,\mu) \mbox{ is recurrent}\}\subset \R$ is open.
In other words for any recurrent BMC with mean offspring $m$ there
exists some $\eps>0$ such that the BMC with mean offspring
$m-\eps$ is still recurrent and hence the critical BMC must be
transient. So assume the BMC to be recurrent. Due to Theorem
\ref{thm:1} \emph{(i)}
 there exists some $k$ such that $\E_\orig
\eta(k,\orig)>1.$ We define an embedded Galton-Watson process
($\zeta_i)_{i\geq 0}.$ We start a BMC* with origin $\orig$ with
one particle in $\orig.$ Let $\Psi_1$ be the particles that are
the first particles in their ancestry line to return to $\orig$
before time $k.$ We define $\Psi_i$ inductively as the number of
particles that have an ancestor in $\Psi_{i-1}$ and are the first
in the ancestry line of this ancestor to return to $\orig$ in at
most $k$ time steps. Clearly $\zeta_0:=1$ and
$\zeta_i:=|\Psi_i|,~i\geq 1,$ defines a supercritical
Galton-Watson process since $E\zeta_1=\E_\orig\eta(k,\orig)>1.$
Furthermore, $E\zeta_1=\E_\orig \eta(k,\orig)$ and
$$ \E_x \eta(k,\orig)= m\cdot \sum_{y\neq \orig} p(x,y) \E_y \eta(k-1,\orig) +
m\cdot p(x,\orig).$$ Now, it is easy to see that  $E\zeta_1$ is
continuous in $m.$ Eventually, for $\eps>0$ sufficiently small the
BMC with mean offspring $m-\eps$ is still recurrent.

\end{proof}

\begin{rem}\label{rem:green}
Theorem \ref{thm:1a} implies that $G(x,x|m)<\infty$ is equivalent
to  transience of the process if and only if
$G(x,x|1/\rho(P))<\infty,$ i.e., the underlying Markov chain is
$\rho$-transient.
\end{rem}

\begin{rem}\label{rem:seed}
The fact that $m \rho(P)>1$ implies the recurrence of the BMC can
be also seen by dint of the interpretation as a general branching
process and Lemma \ref{specapprox}. Due to the latter there exists
a finite and irreducible $Y$ such that $m \rho(P_Y)>1.$ Now, let
us consider only particles in $Y$ and neglect the particles
leaving $Y.$ This defines a supercritical multi-type Galton-Watson
process with first moments $m\cdot P_Y$ that survives with
positive probability, compare with  Chapter V in \cite{athreya},
and hence $\alpha(x)>0.$ Finite regions as $Y,$ that may produce
an infinite number of particles without help from outside, are
referred to as \emph{seeds.}  Note that in \cite{comets:05}
regions  of these kind are called  \emph{recurrent seeds.}
\end{rem}
\subsection{Examples}\label{subsec:ex}

\begin{ex}\label{ex:BRWonZ^d}
Consider the random walk on $\Z^d,$ $d\in\N.$
Let  $e_i\in\Z^d$ with $(e_{i})_j=\delta_{ij}$ for
$i,j\in\{1,\ldots,d\},~ d\geq 1,$ and define transition probabilities  $P$  by
$$ p(x,x+e_i)=p_i^+,~ p(x,x-e_i)=p_i^- ~\mbox{ such that } $$ $$\sum_{i=1}^d p_i^+
 +\sum_{i=1}^d p_i^-=1,\quad\forall x\in\Z^d$$
and such that $P$ is irreducible. Take branching distributions
with constant mean offspring $m$. We calculate, using for example
large deviation estimates (compare with Lemma
\ref{def:spectr_rad}):
$$\rho(P)=2\sum_{i=1}^d \sqrt{p_i^+ p_i^-}.$$
Hence, the corresponding BMC is  transient if and only if
    $$ m \leq 1/\left({2\sum_{i=1}^d \sqrt{p_i^+ p_i^-}}\right)\, .$$
In particular, the BRW $(\Z^d,\mu)$ with constant mean offspring is recurrent if $m>1.$
\end{ex}

\begin{ex}\label{ex:amen}
We consider an irreducible symmetric random walk on a finitely
generated group and constant mean offspring $m.$ We  can classify
groups in amenable or nonamenable using branching random walks: a
finitely generated group $G$ is amenable if and only if every BRW
on $G$  with constant mean offspring $m>1$ is recurrent. This
statement is a variation  of Proposition 1.5 in
\cite{benjamini:94b} where tree-indexed Markov chains are
considered. To proof it, we merely need to combine Theorem
\ref{thm:1a} with the well-known result of Kesten stating that
every irreducible and symmetric Random Walk on a finitely
generated group $G$ has spectral radius $1$ if and only if $G$ is
amenable, e.g. compare with Corollary 12.5 in \cite{woess}.
\end{ex}

\begin{ex}\label{ex:BRWontree}
Let $(X,P)$ be the SRW on the regular tree $\T_M.$ We have
$\rho(P)=\rho(\T_M)=\frac{2\sqrt{M-1}}M$ (compare with  Lemma 1.24
in \cite{woess}). The BMC $(\T_M,\mu)$ with constant offspring
distribution $\mu$ is transient if and only if
$$ m\leq \frac{M}{2\sqrt{M-1}}.$$
\end{ex}

\begin{ex}\label{ex:counter1}
We consider the example of Section 5 in  \cite{comets98} on $\Z$
with binary branching, i.e., $\mu_2(x)=1$ for all $x\in\Z.$ The
transition probabilities are $p(1,0)=p(1,2)=1/8,~p(1,1)=3/4$ and
$$p(x,x+1)=1-p(x,x-1)=\frac{2+\sqrt{3}}4\quad x\neq 1.$$

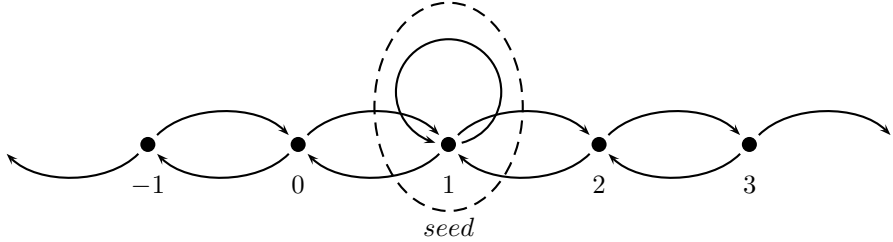
\begin{center}
\begin{pspicture}(-3,-1.5)(7,2)

\multirput*(-2,0)(2,0){5}{\pscircle*{0.1}}
%\psline[linewidth=0.05]{-}(-2.5,0)(2.5,0)
\psset{nodesep=5pt} \psset{arcangle=45}
 \pnode(0,0){0} \pnode(2,0){1} \pnode(6,0){3}
\pnode(4,0){2}\pnode(-2,0){-1} \pnode(8,0){4}\pnode(-4,0){-2}

\ncarc{->}{0}{1} \ncarc{->}{1}{2} \ncarc{->}{1}{0}
\ncarc{->}{2}{1} \ncarc{->}{-1}{0}\ncarc{->}{0}{-1}
\ncarc{->}{2}{3}\ncarc{->}{3}{2}
\ncarc{->}{3}{4}\ncarc{->}{-1}{-2}

\nccircle{->}{1}{0.7}
%\ncline[nodesep=0.15]{->}{A}{B}

%\uput[u](0.5,0.5){$\frac12$} \uput[u](-.5,0.5){$\frac12$}

%\ncline[nodesep=0.15]{->}{A}{C}
\psellipse[linestyle=dashed](2,0.5)(1,1.4)

\uput[d](0,-0.25){$0$} \uput[d](2,-0.25){$1$}
\uput[d](-2,-0.25){$-1$} \uput[d](4,-0.25){$2$}
\uput[d](6,-0.25){$3$} \uput[d](2,-0.8){$seed$}
\end{pspicture}
\end{center}

\noindent The BMC $(X,P,\mu)$ is not strongly recurrent since the
spatially homogeneous BMC with $m(x)=2$ and
$p(x,x+1)=1-p(x,x-1)=\frac{2+\sqrt{3}}4$ for all $ x\in \Z$ is
transient, see Example \ref{ex:BRWonZ^d}. Let us first take $0$ as
the origin $\orig$ of the corresponding BMC*. We show that
$E_1\nu(0)=\infty.$ The number of particles which never leave
state $1$ is given by a Galton-Watson process with mean number
$2\cdot 3/4>1.$ And so, with positive probability, an infinite
number of particles visits  state $1.$ This clearly implies
$\E_1\nu(0)=\infty.$ Eventually, the BMC $(X,P,\mu)$ is recurrent,
but not strongly recurrent. Notice, if  $\orig=1$ then we have for
the corresponding BMC* that $\E_x\nu(1)<\infty$ for all $x.$
\end{ex}

\begin{rem}\label{rem:tran_local}
Example \ref{ex:counter1} illustrates very well the idea of
\emph{seeds} that make BMCs recurrent:   state $1$ can be seen as
a seed that may create an infinite number of particles without
help from outside. Since in this case the seed is just a
\emph{local inhomogeneity,} the process may \emph{escape} the seed
and  is not strongly recurrent.
\end{rem}

\section{Criteria for strong recurrence}\label{sec:strong}
In this section we  discuss criteria for strong recurrence of BMC.
In Subsection \ref{sec:gen_crit} we present known and new
conditions for general Markov chains and offspring distributions.
While  Theorem \ref{thm:strecur} of \cite{menshikov:97} is more of
theoretical interest the  Lemma \ref{lem:strongrec} and
Proposition \ref{prop:strong:trivial} are new useful tools to
prove strong recurrence. In particular, if the underlying
 Markov chain and the offspring distributions are \emph{homogeneous,} recurrence
and strong recurrence coincide, compare with Subsection
\ref{subsec:homo}.

In Subsection \ref{subsect:finer} we present several approaches in
order to develop sufficient and necessary conditions for strong
recurrence for general BMC, see Lemma \ref{lem:rec_common_root},
Theorem \ref{thm:rec_quasi}, and Lemma \ref{lem:tilderho}. An
interesting observation is that transience / recurrence depend on
local properties and recurrence / strong recurrence on global
properties of the process. Therefore a classification result would
demand a suitable description of infinite structures and would
deliver a measure for inhomogeneity of the process. The conditions
for strong recurrence are given in terms of appropriate spectral
radii. While a general and applicable criterion for strong
recurrence remains challenging, our condition work well in
concrete situations, e.g. compare with Theorem
\ref{thm:finite_cone} and  Example \ref{ex:line_tree}. The section
ends with a short discussion including conjectures in Subsection
\ref{sec:outlook}.

\subsection{General Criteria}\label{sec:gen_crit}
The criteria for transience and recurrence, see Theorems
\ref{thm:1} and \ref{thm:1a}, do not depend on the precise
structure of the branching mechanism but only on the mean
offspring $m.$  This can no longer hold true for criteria for
strong recurrence since we can choose the branching distribution
such that with positive probability no branching occurs at all,
see the following Example \ref{ex:gen_crit}.

\begin{ex}\label{ex:gen_crit}
Consider the random walk on $\Z$ with drift to the right, i.e.,
$p(x,x+1)=1-p(x,x-1)=p>1/2~\forall x\in\Z.$ In order to construct a nontrivial
BMC where with positive probability no branching occurs, i.e.,
$\P(\eta(n)=1~\forall n\geq 1)>0,$ we first investigate the
underlying random walk. We know, Law of Large Numbers, that $S_n/n\to s:=2p-1$ as $n\to\infty.$ Hence for
each realization $S_n(\omega)$ of the random walk there exists
some $T(\omega)$ such that $S_n(\omega)>(s-\eps)n$ for all $n>
T(\omega)$ for some small $\eps>0.$  Define
$$C_T:=\{\omega: S_n(\omega)>(s-\eps)n~\forall n>T\} \mbox{ and
} C_\infty:=\bigcup_{T=1}^\infty C_T.$$ Due to the Law of Large
Numbers and since the $C_T$ are increasing, we have
$1=\P(C_\infty)=\lim_{T\to\infty} \P(C_T).$ Hence, there exists
some $T>0$ such that $\P(A)>0,$ with $A:=\{\omega:
S_n(\omega)>(s-\eps)n~\forall n>T\}.$ Now we choose the branching
distributions such that on the event $A$ with positive probability
no branching occurs. We define $\left(\mu(x)\right)_{x\in\Z}$ such
that $m(x)=m>1/\rho(P)$ and $\mu_1(x)=1-e^{-b x}$ for $x>0,$ and
$\mu_1(x)=(1-e^{-b})$ for $x\leq 0,$ where $b$ is some positive
constant. Eventually,
$$\P(\eta(n)=1~ \forall n\geq 1|A)\geq \left(1-e^{-b}\right)^T
\prod_{n=T}^\infty \left(1-e^{-b (s-\eps)n}\right)>0$$ and the BMC
$(X, P, \mu)$ is not strongly recurrent but recurrent since
$m>1/\rho(P).$ On the other hand if $\tilde\mu(x)=\tilde\mu$ with
$m>1/\rho(P)$  and hence is homogeneous, then  the BMC $(X, P,
\tilde\mu)$ is strongly recurrent, compare with Subsection
\ref{subsec:homo}.
\end{ex}

Despite the above discussion, there exists a sufficient and
necessary condition for strong recurrence where the offspring
distribution may depend on the states. Let
$$\Psi(x,z):=\sum_{k=1}^\infty z^k \mu_k(x)$$ be the generating
function of $\mu(x).$ We have the following necessary and
sufficient condition for strong recurrence of \cite{menshikov:97}.

\begin{thm}\label{thm:strecur}
The BMC $(X,P,\mu)$ is not strongly recurrent if and only if there
exists a finite subset $M$ of $X$ and a function $0<g\leq 1,$ such
that

\begin{equation}\label{eq:thm:strecur:1}
    \Psi\left(x,Pg(x)\right)\geq g(x) \quad \forall x\notin
     M
\end{equation} and
\begin{equation}\label{eq:thm:strecur:2}
    \exists y\notin M:~ g(y)>\max_{x\in M} g(x).
\end{equation}

\end{thm}

\begin{proof}
We give a sketch of the proof in  \cite{menshikov:97}.  We start the BMC in $y$ and define
$$\widetilde Q(n):=\prod_{i=1}^{\eta(n)} g(x_i(n)).$$ Furthermore, let $$\tau:=\min_{n\geq 0} \{\exists
i\in\{1,\ldots,\eta(n)\}:~x_i(n)\in M\}$$
be the entrance time in $M.$   It turns out that
$$Q(n):=\widetilde Q(n\wedge \tau)$$ is a submartingal for $n\geq
0.$  Since $Q(n)$ is bounded, it converges a.s. and in $L^1.$ Hence there
exists some random variable $Q_\infty$ such that
$$Q_\infty=\lim_{n\to\infty} Q(n)$$ and
\begin{eqnarray}\label{eq:thm:strecu}
  \E_{y}Q_\infty &=& \lim_{n\to\infty} \E_{y}Q(n)\geq \E_{y}Q(0)=g(y).
\end{eqnarray}
Assuming that the BMC is strongly recurrent, we obtain that
$\tau<\infty$ a.s. and therefore $Q_\infty\leq \max_{x\in M}
g(x).$ This contradicts inequality (\ref{eq:thm:strecu}) since
$g(y)>\max_{x\in M} g(x).$ The converse is more constructive.
Assume the BMC not to be  strongly recurrent and consider the
probability that starting the BMC in $x$ no particles hits some
$\orig\in X:$
$$g(x):=\P_x\left(\forall n\geq 0~\forall i=1,\ldots,\eta(n):
~x_i(n)\neq \orig \right)\quad \mbox{ for } x\neq \orig$$ and
$g(\orig):=0.$  One easily checks that $g$ verifies the requested conditions for
$M:=\{\orig\}.$
\end{proof}

The conditions in Theorem \ref{thm:strecur} are difficult to
check. We did not find a more explicit formulation. Furthermore, it
is not clear if  strong recurrence depends on the whole structure
of the branching distributions, since the above conditions are
written in term of the generating function of $\mu.$ Nevertheless,
we see in Subsection \ref{subsec:homo} that in homogenous cases
the necessary and sufficient condition for strong recurrence does
only depend on the mean offspring $m(x)$ and conjecture that this
fact holds in general, see Conjecture \ref{conj:strongrec2}.

\begin{rem}
Theorem \ref{thm:strecur} implies in particular that strong
recurrence depends on global properties of the BMC since local
properties can be excluded by the choice of the finite set $M.$
\end{rem}

A useful tool are induced Markov chains that couple a Markov chain
to the branching Markov chain. The induced Markov chain $X_n$ is
defined inductively. We can think of it as the position of a
label. The BMC starts, at time $n=0,$ with one particle that is
labelled. At time $n$ the labelled particle picks at random one of
its offspring and hand over the label. It is easy to check that
the position of the label defines a Markov chain with transition
probabilities $P.$ Another way to interpret the induced Markov
chain is to modify the original process in a way that particles
do not die but produce offspring with distribution
$\tilde\mu_{i-1}=\mu_{i},~i\geq 1.$ In this case we can speak of
the trace of a particle which has the same distribution as the
trace of a Markov chain with transition kernel $P.$

The next Lemma is our main tool to show strong recurrence.
\begin{lem}\label{lem:strongrec} Let $c>0$ and define $C:=\{\alpha(x)\geq
c\}.$ If the set $C$ is recurrent with respect to the Markov chain
$(X,P)$, i.e.,  is a.s. hit infinitely often by the trace of the
Markov chain, the BMC is strongly recurrent.
\end{lem}
\begin{proof}
The idea is to define a sequence of embedded supercritical
Galton-Watson processes  and show that at least one of them
survives. We start the process with $x=x_1\in C.$ Let us define
the first Galton-Watson process $(\zeta_i^{(1)})_{i\geq 0}.$ To
this end, let $\Psi_1$ be the particles that are the first
particles in their ancestry line to return to $x$ before time $k$
(to be chosen later) and  define $\Psi_i$ inductively as the
number of particles that have an ancestor in $\Psi_{i-1}$ and are
the first in the ancestry line of this ancestor to return to $x$
in at most $k$ time steps. Clearly $\zeta_0^{(1)}:=1$ and
$\zeta_i^{(1)}:=|\Psi_i|,~i\geq 1,$ defines a Galton-Watson
process. Due to the definition of the process BMC*
we have that $E\zeta_1^{(1)}=\E_x\eta(k,x).$\\
\noindent {\bf Claim:} There is some $k$ such that
$E\zeta_1^{(1)}>1$ and that the probability of survival of
$(\zeta_i^{(1)})_{i\geq 0}$ is larger than $ c/2.$

We choose $k$ such that the probability of survival of
$(\zeta_i^{(1)})_{i\geq 0}$ is larger than $ c/2.$  If this first
Galton-Watson process dies out we wait until the induced Markov
chain hits a point $x_2\in C$; this happens with probability one
since $C$ is recurrent w.r.t the induced Markov chain. Then we
start a second process $(\zeta_i^{(2)})_{i\geq 0},$ defined in the
same way as the first but started in  position $x_2.$ If the
second process dies out, we construct a third one, and so on. We
obtain a sequence of independent Galton-Watson processes
$\left((\zeta_i^{(j)})_{i\geq 0}\right)_{j\geq 1}.$ The
probability that all these processes die out is less than
$\prod_{j=1}^\infty (1-c/2)=0.$ Eventually, at least one process
survives and we have $\alpha(x)=1$ for all $x\in X.$

It remains to prove the claim. Consider the  Galton-Watson process
$(Z_i)_{i\geq 0}$ constructed as $(\zeta_i^{(1)})_{i\geq 0}$ but
with $k=\infty.$ Hence $EZ_1=\E_x \nu(x)\in (0,\infty].$ Let
$f(s)=\sum_{j=1}^\infty \mu_j s^j,~|s|\leq 1$ be the generation
function of $(Z_i)_{i\geq 0}.$ From the definition of $f$ as a
power series with nonnegative coefficients, we have that it is
strictly convex and increasing in $[0,1).$ Furthermore, the
extinction probability $q$ of $(Z_i)_{i\geq 0}$ is the smallest
nonnegative root of the equation $t=f(t).$ For every $k$ we define
a process $(Z^k_i)_{i\geq 0}$ with corresponding mean offspring
$\eta(k,x),$ distribution $\mu^k=(\mu_1^k,\mu_2^k,\ldots),$ and
generating function $f^k.$ The probabilities $\mu^k$ converge
pointwise to $\mu$ and so do the generating functions. Using the
fact that $f(q)=q$ and $1-q\geq c$ we find a $k$ such that
$f^k(1-c/2)\leq 1- c/2,$ thus $(Z^k_i)_{i\geq 0}$ survives with
probability at least $c/2.$
\end{proof}

\begin{rem}\label{rem:lem:strongrec}
In Lemma \ref{lem:strongrec} we can replace the condition that $C$
is recurrent w.r.t. the Markov chain by the condition that $C$
is recurrent w.r.t. the BMC, i.e., $C$ is visited infinitely
often by some particles of the BMC.
\end{rem}

Let $F(x)$ denote the return probability of the Markov chain $(X,
P)$, i.e., the probability that the Markov chain  started in $x$
returns to $x.$ If we assume the branching distributions to be
constant, i.e., $\mu(x)=\mu$ for all $x\in X,$ we have  the
following sufficient condition for  strong recurrence in terms of
the mean offspring and the return probability of the underlying
Markov chain.
\begin{prop}\label{prop:strong:trivial}
The BMC $(X,P,\mu)$ with constant offspring distribution
%, i.e., $\mu(x)=\mu$ for all $x\in X,$ 
is strongly recurrent if
$$m>\sup_{x\in X} \frac1{F(x)}.$$
\end{prop}
\begin{proof}
Due to Lemma \ref{lem:strongrec} we have to show that
$\alpha(x)\geq c$ for all $x\in X$ and some $c>0.$ Consider  the
Galton-Watson process $(\tilde\xi_i)_{i\geq 0}$ with offspring
distribution $\mu$ and mean $m.$ Furthermore, let $p$ such that
$$\frac1m<p<\inf_{x\in X} F(x)$$
and percolate the process $(\tilde\xi_i)_{i\geq 0}$ with survival
parameter $p.$ This leads to a Galton-Watson process
$(\xi_i)_{i\geq 0}$ with mean $mp>1$ and some survival probability
$c>0,$ compare with \cite{lyons90}. Back on BMC, we start the
process $(X,P,\mu)$ with one particle in some arbitrary position,
say $\orig,$ and compare the original process with the BMC
$(X,P,\tilde\mu)$ with fewer particles:
$\tilde\mu(\orig):=\mu(\orig)$ and $\tilde\mu_1(x):=1$ for all
$x\neq\orig.$ In other words, $(X, P, \tilde\mu)$ does only branch
in $\orig.$ Observe that the number of particles returning to
$\orig$ in this process can be described by dint of a percolated
Galton-Watson process $(\zeta_i)_{i\geq 0}$ with offspring
distribution $\mu$ and survival parameter $F(\orig).$ Since
$F(\orig)<p$  we can use a standard coupling of Bernoulli
percolation, compare with Chapter 4 in \cite{lyons:book}, to prove
that the survival probability of $(\zeta_i)_{i\geq 0}$ is at least
the one of $(\xi_i)_{i\geq 0}.$ If $(\zeta_i)_{i\geq 0}$ survives,
an infinite number of particles visits $\orig$ in $(X,
P,\tilde\mu)$ and hence in $(X, P,\mu)$ as well. We can conclude
that $\alpha(\orig)\geq c$ for the original BMC $(X, P, \mu).$
\end{proof}

\subsection{Homogeneous BMC}\label{subsec:homo}
Lemma \ref{lem:strongrec} offers a general argument to show strong
recurrence.  In particular,  it is used to prove that homogeneous
BMC are strong recurrent if and only if they are recurrent. This
fact is also plausible from the viewpoint of seeds. An infinite
number of seeds are visited and each of these gives birth to a
supercritical multi-type Galton-Watson process with extinction
probability bounded from below. We give several known
(\ref{sec:quasi}, \ref{sec:brwre}) and new (\ref{sec:conetye},
\ref{sec:percolation}, and \ref{sec:uniform}) examples of
homogeneous processes. They are all consequences of Theorem
\ref{thm:1a} and Lemma \ref{lem:strongrec}.

\subsubsection{Quasi-transitive BMC}\label{sec:quasi}
Let $X$ be a locally finite, connected graph with discrete metric
$d.$ An  automorphism of $X$ is a self-isometry of $X$ with
respect to $d,$ and $AUT(X)$ is the group of all automorphisms of
$X.$ Recall  that when a group $\Gamma$ acts on a set $X,$ this
process is called a group action: it permutes the elements of $X.$
The group orbit of an element $x$ is defined as $\Gamma x :=
\{\gamma x:~ \gamma\in \Gamma\}.$ A group $\Gamma$ acts
transitivly on $X$ if it possesses only a single group orbit,
i.e., for every pair of elements $x$ and $y$ of $X$, there is a
group element $\gamma\in \Gamma$ such that $\gamma x = y.$ The
graph $X$ is called  transitive  if $AUT(X)$ acts transitively on
$X$, and quasi-transitive if $AUT(X)$ acts with a finite number of
orbits. Let $P$ be the transition matrix of an irreducible  random
walk on $X$ and $AUT(X,P)$ be the group of all $\gamma\in AUT(X)$
which satisfy $p(\gamma x,\gamma y)=p(x,y)$ for all $x,y\in X.$ We
say the Markov chain $(X,P)$ is transitive if the group $AUT(X,P)$
acts transitively on $X$ and quasi-transitive if $AUT(X,P)$ acts
with a finite number of orbits on $X.$

The definition of quasi-transitivity can be extended to BMC. We
say a BMC is quasi-transitive if the group $AUT(X,P,\mu)$ of all
$\gamma\in AUT(X,P)$ which satisfy $\mu_k(x)=\mu_k(\gamma
x)~\forall k\geq 1$ for all $x\in X$ acts with a finite number of
orbits on $X.$ Observing that $\alpha(x)$ attains only a finite
number of values and hence $\alpha(x)\geq c$ for some $c>0$ we
obtain due to Theorem \ref{thm:1a} and Lemma \ref{lem:strongrec}
the following result. It is due to \cite{gantert:04} and also
generalizes some results of \cite{stacey:03}.

\begin{thm}\label{thm:2}
Let  $(X,P,\mu)$ be a  quasi-transitive BMC with constant mean
offspring $m(x)= m  >1.$ It holds that
\begin{itemize}
    \item the BMC is transient $(\alpha=0)$ if $m\leq 1/\rho(P)$.
    \item the BMC is strongly recurrent $(\alpha=1)$ if $m>1/\rho(P)$.
\end{itemize}
\end{thm}

\subsubsection{Branching random walk on trees with finitely many cone
types}\label{sec:conetye} An important class of \emph{homogeneous}
trees are  periodic trees that are also known as trees with
finitely many cone types, compare with \cite{lyons:book} and
\cite{nagnibeda:02} . These trees arise as the directed cover
(based on $r$) of some finite connected directed graph $G:$ the
tree $T$ has as vertices the finite paths in $G,$ i.e., $\langle
r,i_1,\ldots,i_n\rangle.$ We join two vertices in $T$ by an edge
if one path is an extension by one vertex of the other. The cone
$T^x$ at $x\in T$ is the subtree of $T$ rooted at $x$ and spannend
by all vertices $y$ such that $x$ lies on the geodesic from $r$ to
$y$. We say that $T^x$ and $T^y$ have the same cone type if they
are isomorphic as rooted trees and every cone type corresponds in
a natural way to some vertex in $G.$ Let $\tau(x)$ be the function
that maps a vertex $x\in T$ to its cone type in $G.$ If $G$ is
strongly connected, i.e., for every pair $x,y$ there is a directed
path in $G$ from $x$ to $y$, we call the cone types irreducible.
In this case every cone \emph{contains} every other cone type as a
subtree.

We consider the nearest neighbour random walk on $T$ according to
\cite{nagnibeda:02}. Suppose we are given transition probabilities
$q(i,j)_{i,j\in G}$ on $G.$ We may hereby assume, w.l.o.g., that
$q(i,j)>0$ if and only if there is an edge from $i$ to $j$ in $G$.
Furthermore, suppose we are given \textit{backward probabilities}
$p(-i)\in (0,1)$ for each $i\in G$. Then the random walk on the
tree $T$ is defined through the following transition probabilities
$p(x,y)$, where $x,y\in T$:
\begin{equation*}
p(\orig,y):= q(r,\tau(y)), \mbox{ if } x=y^-,
\end{equation*}
and for $x\neq \orig$ with $\tau(x)=i$
\begin{equation*}
p(x,y):=
\begin{cases}
\bigl(1-p(-i)\bigr) q\bigl(\tau(x),\tau(y)\bigr), & \mbox{ if }x=y^-  \\
p(-i), &  \mbox{ if }y=x^-,
\end{cases}
\end{equation*}
where $x^-$ is the \emph{ancestor} of $x.$  It's worth to mention,
that  if the cone type are irreducible, then $\rho(P)<1$ if and
only if the random walk is transient, compare with Theorem B in
\cite{nagnibeda:02}.

Furthermore, we assign branching distributions $\mu$ to the
vertices of $G$ and define the BRW $(T,P,\mu)$ with
$\mu(x)=\mu(\tau(x)).$ We speak of a BRW of finitely many
(irreducible) cone types and have the following classification:

\begin{thm}\label{thm:2_finitecone}
Let  $(T,P,\mu)$ be BRW with finitely many irreducible cone types
and constant mean offspring $m(x)= m  >1.$ It holds that the BMC
\begin{itemize}
    \item  is transient $(\alpha=0)$ if $m\leq 1/\rho(P),$
    \item  is strongly recurrent $(\alpha=1)$ if $m>1/\rho(P).$
\end{itemize}
\end{thm}

\begin{proof}
First observe that the process is not quasi-transitive and
$\alpha(x)\neq \alpha(y)$ for $\tau(x)=\tau(y).$ We prove that for
every cone type  $\alpha(x)\geq c(\tau(x))>0$ and conclude with
Theorem \ref{thm:1a} and Lemma \ref{lem:strongrec}. Due to Lemma
\ref{specapprox} there exists a finite subset $F$ such that
$m\cdot \rho(P_F)>1.$ Now, since the cone types are irreducible
every cone $T^x$ contains $F$ as a subset (hereby we mean a subset
that is isomorphic to $F$). We fix a cone type say $i\in G$ and
and let $x$ such that $\tau(x)=i.$ There exists some $n\in\N$ such
that the ball $B_n(x)$ of radius $n$ around $x$ contains $F.$
Lemma \ref{specapprox} yields that $\rho(T^x\cap B_n(x))\cdot
m>1.$ Recalling Remark \ref{rem:seed} we notice that the embedded
multi-type Galton-Watson process living on $T^x\cap B_n(x)$ is
supercritical and survives mit positive probability, say $c(i).$
We can conclude with a standard coupling arguments that
$\alpha(x)\geq c(i)$ for all $\tau(x)=i.$
\end{proof}

If the cone types are not irreducible all  three phases may occur,
compare with Theorem \ref{thm:finite_cone}.

\subsubsection{Branching random walk in random environment (BRWRE) on Cayley
graphs}\label{sec:brwre}

Let $G$ be  a finitely generated group. Unless $G$ is abelian, we
write the group operation multiplicatively. Let $S$ be a finite
symmetric generating set of $G,$ i.e., every element of $G$ can be
expressed as the product of finitely many elements of $S$ and
$S=S^{-1}$. The Cayley graph $X(G,S)$ with respect to $S$ has
vertex set $G$, and two vertices $x,y\in G$ are joined by an edge
if and only if $x^{-1}y\in S.$ Now, let $q$ be some probability
measure on $S.$  The random walk on $X(G,S)$ with transition
probabilities $q$ is the Markov chain with state space $X=G$ and
transition probabilities
$$p(x,y)=q(x^{-1}y)\quad \mbox{for }x^{-1}y\in S$$ and $0$ otherwise.
The discrete convolution is defined as $q*q(x)=\sum_y q(y)
q(y^{-1} x).$ The $n$-step transition probabilities are
$$p^{(n)}(x,y)=q^n(x^{-1}y),$$
where $q^n$ is the $n$-fold discrete convolution of $q$ with
itself. We start the random walk at time $0$ in some position
$\orig\in X.$

We introduce the random environment. Let $\mathcal{M}$ be the
collection of all probability measures on $S$ and let
$(\omega_x)_{x\in X}$ be a collection of iid random variables with
values in $\mathcal{M}$ which serve as an environment. For each
realization $\omega:=(\omega_x)_{x\in X}$ of this environment, we
define a Markov chain $(X_n)_{n\in \N}$ on $X=G$ with starting
position $\orig$ and
$$\P_{\omega,\orig}(X_{n+1}=y| X_n=x):=p_\omega(x,y):=\omega_x(x^{-1} y)\quad
\forall n\geq 1.$$ We denote by $P_\omega$  the transition kernel
of the Markov chain on the state space $X.$

Let $\eta$ be the distribution of this environment. We assume that
$\eta$ is a product measure with  one-dimensional marginal $Q.$
The support of $Q$ is denoted by $\mathcal{K}$ and its convex hull
by $\hat{\mathcal{K}}.$ We always assume the following condition
on $Q$ that ensures the irreducibility of a random walk with
transition probabilities $q\in \hat{\mathcal{K}}:$
\begin{eqnarray}\label{eq:cond}
      Q\{\omega: \omega(s)>\gamma~\forall
    s\in S'\}=1 \mbox{ for some } \gamma>0,
\end{eqnarray} where $S'\subseteq S$ is a minimal set of
generators, i.e., every proper subset $T\subsetneq S'$ is not a
generating set.

In addition to the  environment which determines the random walk
we introduce a random environment determining the branching
mechanism. Let $\mathcal{B}$ be the set of all infinite positive
sequences $\mu=\left(\mu_k\right)_{k\geq 1}$ satisfying
$\sum_{k=1}^\infty \mu_k=1$ and $m(\mu):=\sum_{k=1}^\infty
k\mu_k<\infty.$ Let $\widetilde Q$ be a probability distribution
on $\mathcal{B}$ and set
\begin{equation}
\label{mstar} m^*:=\sup\{m(\mu):\mu\in {\rm supp}(\widetilde Q)\}
 \end{equation}
which may take values in $\R\cup\{\infty\}.$ Let $(\omega_x)_{x\in
X}$ be a collection of iid random variables with values in
$\mathcal{M}$ and $(\mu_x)_{x\in X}$ be a collection of iid random
variables with values in $\mathcal{B}$ such that $(\omega_x)_{x\in
X}$ and $(\mu_x)_{x\in X}$ are independent, too. Let $\Theta$ be
the corresponding product measure with one-dimensional marginal
$Q\times \widetilde Q.$ For each realization
$(\omega,\mu):=(\omega_x,\mu_x)_{x\in X}$ let $P_\omega$ be the
transition kernel of the underlying Markov chain and branching
distribution $\mu(x)=\mu_x.$  Thus, each realization
$(\omega,\mu)$ defines a BMC $(X,P_\omega,\mu).$ We denote by
$\P_{\omega,\mu}$ the corresponding probability measure.

We assume that $m^*>1,$ excluding the case where the BMC is
reduced to a Markov chain without branching.

The classification is due to \cite{mueller06} where it is proved
for BRWRE on Cayley graphs. Furthermore, compare to
\cite{mueller07} where it is shown for a model where branching and
movement may be dependent. The interesting fact is that the type
only depends on some extremal points of the support of the random
environment, namely the highest mean offspring and the \emph{less
transient} homogeneous random walk.

We obtain due to Lemma \ref{lem:strongrec} and Theorem
\ref{thm:1a} that the spectral radius is deterministic, i.e.,
$\rho(P_\omega)=\rho$ for $\Theta$-a.a. realizations. This can
also be seen directly from the observation that
$\rho(P_\omega)=\limsup \left( p^{(n)}(x,x)\right)^{1/n}$ does not
depend on $x$ and hence by ergodicity of the environment is
constant a.s. .

\begin{thm}\label{thm:5}
 If $m^*\leq  1/{\rho} $ then the BRWRE is transient for
$\Theta$-a.a. realizations $(\omega,\mu)$, otherwise it is
strongly recurrent for $\Theta$-a.a. realizations $(\omega,\mu)$.
\end{thm}

In the special case of the lattice the spectral radius $\rho$ can
be calculated explicitly.

\begin{cor}\label{cor:6}
The BRWRE on $\Z^d$ is strongly recurrent for $\Theta$-a.a. realizations if
$$ (m^*)^{-1} <\sup_{p\in\hat{\mathcal{K}}} \inf_{\theta\in\R^d} \left( \sum_s e^{\langle
   \theta, s\rangle} p(s)\right).$$ Otherwise it is transient for
   $\Theta$-a.a. realizations.
\end{cor}

\subsubsection{BRW on percolation clusters}\label{sec:percolation}

Let in this subsection $G$ be a  graph with bounded geometry and
origin $\orig.$ We consider Bernoulli$(p)$ percolation on $G,$
i.e., for fixed $p\in [0,1],$ each edge is kept with probability
$p$ and removed otherwise, independently of the other edges.
Denote the random subgraph of $G$ that remains by $C(\omega)$ and
$C(\omega,x)$ the connected component containing $x.$ We refer  to
Chapter $6$ in \cite{lyons:book} for more information and
references on percolation models.

\begin{thm}\label{thm:7}
The BRW with constant offspring distribution and $m>1$ and
underlying SRW on a connected component of $C(\omega)$ is a.s.
strongly recurrent.
\end{thm}

\begin{proof}
We start the BRW in $\orig$ and consider $C(\omega,\orig).$
Clearly if the component is finite then the BRW is strongly
recurrent. Due to Lemma \ref{specapprox} and  $\rho(\Z)=1$ there
exists a subset $Y$ of $\Z$ such that $\rho(P_Y)\cdot m >1.$
W.l.o.g. we can assume $Y$ to be a line segment of length $k.$
Now, let us imagine that the percolation is constructed during the
evolution of the BRW. For $n\ge 1$ we denote $B_n=B_n(\orig)$ the ball of radius $n$ around the origin $\orig.$ We percolate the edges in the ball
$B_k.$ The percolation
cluster is now defined inductively. If one vertex of the border,
say $x_i,$ of the ball $B_{ik}$ is hit by some particle we percolate
the edges in $B_{(i+1)k}\setminus B_{ik}.$ With positive
probability $\beta(x_i)$ we have that $C(\omega,x_i)\cap
(B_{(i+1)k}\setminus B_{ik})$ equals the line segment of length
$k.$ Since $G$ is of bounded geometry we have that
$\beta(x_i)\geq\delta>0$ for all $x_i.$ Observe that $\alpha(x_i)$
is at least the survival probability, say $c,$ of the multi-type
Galton-Watson process restricted on the line segment of length
$k,$ compare with Remark \ref{rem:seed}. Eventually, either
$C(\omega,\orig)$ is finite or the set $\{\alpha(x)\geq
c\cdot\delta\}$ is recurrent w.r.t. the BRW and we can conclude
with Remark \ref{rem:lem:strongrec}.
\end{proof}

\subsubsection{Uniform BMC}\label{sec:uniform}
Let us assume that $(p^{(l)}(x,x))^{1/l}$ converges uniformly in
$x$, i.e., $\forall~\eps>0~\exists~l$ such that $
(p^{(l)}(x,x))^{1/l}>\rho(P)-\eps~\forall x\in X,$
 and that there is a $k\in\N$ such that
$\inf_x \sum_{i=1}^k i \mu_i(x) \geq 1/\rho(P).$ Now,  consider a
modified  BMC with branching distributions
$$\tilde\mu_0(x)=\sum_{i=k+1}^\infty \mu_i(x)\mbox{ and }
\tilde\mu_i(x)=\mu_i(x)\mbox{ for }i=1,\ldots,k\mbox{ and }x\in
X.$$ For this new process we obtain a sequence of supercritical
Galton-Watson processes $((\zeta_i^{(j)})_{i\geq 1})_{j\geq 1}$
with bounded variances and means bounded away from $1$, since $l$
and $k$ do not depend on the starting position ${s_j}.$ Observe
that we have for  the generating function $f$ of a Galton Watson
process with mean $m$ and variance $\sigma^2$ that $f'(1)=m$ and
$f''(1)=\sigma^2/(m-m^2).$ The extinction probability $q$ of a
Galton-Watson process is the unique nonnegative solution less than
$1$ of the equation $s=f(s).$ Using Taylor's Theorem and the
convexity of $f'$ we can conclude that the extinction
probabilities $q_j$ of $(\zeta_i^{(j)})_{i\geq 1}$ are bounded
away from $1.$

\subsection{Inhomogeneous BMC}\label{subsect:finer}

In this section we give conditions for strong recurrence,
$\alpha\equiv 1,$ and recurrence, $\alpha<1,$  that, although
failing to produce a complete classification,  work well in
concrete examples. We assume throughout this section that
$\mu(x)=\mu$ for all $x\in X.$

\subsubsection{Connecting Markov chains at a common
root}\label{subsec:glue} We present a method to glue different
Markov chains following Chapter 9 in \cite{woess}.  Let
$(X_i,P_i), i\in I,$ be a family of irreducible Markov chains. We
choose a root $r_i$ in each $X_i$ and connect the $X_i$ by
identifying all these roots. The rest of the $X_i$ remains
disjoint. This gives a set $X=\bigcup_i X_i$ with root
$r,~\{r\}=\bigcap_i X_i.$ In order to define the transition matrix
$P$ on $X$, we choose constants $\alpha_i>0$ such that $\sum_i
\alpha_i=1$ and set
\begin{eqnarray}\label{eq:glueMC}
  p(x,y) &=& \left\{ \begin{array}{ll}
    p_i(x,y) & x,y\in X_i,~x\neq r, \\
    \alpha_i p_i(r,y) & x=r,~ y\in X_i\setminus\{r\}, \\
    \sum_i \alpha_i p_i(r,r) & x=y=r,\\
    0 & \mbox{otherwise.}
  \end{array}\right.
\end{eqnarray}

When each $X_i$ is a graph and $P_i$ is the SRW on $X_i,$ then $X$
is the graph obtained by connecting the $X_i$ at a common root
$r.$ Choosing $\alpha_i=deg_{X_i}(r)/deg_X(r),$ we obtain the SRW
on $X.$ Due to this construction the calculation of the spectral
radius of $(X,P)$ can be done with the help of  generating
functions of $(X_i,P_i),$ compare with Chapter 9 in \cite{woess}.

For these types of Markov chains we obtain a condition for strong
recurrence in terms of
\begin{equation}\label{eq:vr}
\vr(P):=\inf_{i\in I} \rho(P_i)\in[0,1].
\end{equation}
\begin{lem}\label{lem:rec_common_root}
Let $(X_i,P_i), i\in I,$ be a family of irreducible Markov chains
and $(X,P)$ as defined in (\ref{eq:glueMC}). The BMC $(X,P,\mu)$
with constant branching distribution is not strongly recurrent, i.e., $\alpha<1,$ if
$$m< 1/\vr(P).$$ If the $\inf$ is attained then $m=1/\vr(P)$ implies
$\alpha<1,$ too.
\end{lem}
\begin{proof}
There exists $i\in I$ such that $m\leq 1/\rho(P_i).$  Due to
Theorem \ref{thm:1a} we know that the BMC $(X_i,P_i,\mu)$  is
transient. Hence,  there exists some $x\in X_i$ such that the BMC
$(X_i,P_i,\mu)$ started in $x$ never hits $r_i=r$ with positive
probability.  Therefore, with positive probability the BMC
$(X,P,\mu)$ started in $x$ never hits $r.$
\end{proof}

\begin{rem}
If $m=1/\vr(P)$ and the  $\inf_{i\in I} \rho(P_i)$ is not attained
then both cases can occur. The BMC is strongly recurrent if all
$(X_i,P_i)$ are quasi-transitive, compare with Theorem
\ref{thm:rec_quasi}. In order to construct an example that is not
strongly recurrent, let $(X_1,P_1)$ be as in Example
\ref{ex:counter1}.  For $i\geq 2,$ let $(X_i,P_i)$ be the random
walk on $\Z$ with drift defined by
$p_i(x,x+1)=\frac{2+\sqrt{3}}4-\frac1{i+1}.$ We glue the Markov
chains in $r=r_i=0$ and obtain $\vr(P)=\frac12.$  Since the BMC
$(X_1,P_1,\mu)$ with $m=2$ is not strongly recurrent, this follows
for the BMC $(X, P, \mu)$ as well.
\end{rem}

For certain Markov chains, constructed as above, we can give a
complete classification in  transience, recurrence and strong
recurrence. Observe that we can replace \emph{quasi-transitive} by
any other homogeneous process of Subsection \ref{subsec:homo}.
Interesting is the subtle behavior in the second critical value;
the BMC may be strongly recurrent or weakly recurrent.

\begin{thm}\label{thm:rec_quasi}
Let $(X_i,P_i), i\in I,$ be a family of quasi-transitive
irreducible Markov chains and $(X,P)$ as defined in
(\ref{eq:glueMC}). We have the following classification for the
BMC $(X,P,\mu)$ with constant mean offspring $m:$
\begin{equation}
  \begin{array}{clcrcl}
  (i)& &m& \leq 1/\rho(P) &\Longleftrightarrow & \alpha\equiv 0, \\
  (ii)& 1/\rho(P) < &m& < 1/\vr(P)    &\Longleftrightarrow & 0< \alpha(x) < 1,\\
  (iii)& 1/\vr(P) < &m& &\Longleftrightarrow & \alpha\equiv 1.
  \end{array}
\end{equation}

If the $\inf$ in the definition of $\vr(P)$ is attained, then
$m=1/\vr(P)$ implies $\alpha<1,$ and if the $\inf$ is not
attained, then $m=1/\vr(P)$ implies that $\alpha\equiv 1.$
\end{thm}
\begin{proof}
The part $(i)$ is Theorem \ref{thm:1a},  $(ii)$ is Lemma
\ref{lem:rec_common_root} and $(iii)$ follows from Theorem
\ref{thm:2} by observing that each BMC $(X_i,P_i,\mu)$ is strongly
recurrent. The same argumentation holds if the $\inf$ is not
attained. The case when the $\inf$ is attained follows with Lemma
\ref{lem:rec_common_root}.
\end{proof}

Analogous arguments yield the classification for trees with
finitely many cone types that are not necessarily irreducible. For
this purpose let $G_i$ be the irreducible classes of $G,$ $T_i$
the directed cover of $G_i,$ and $\widetilde\rho(T):=\min_i
\rho(T_i).$

\begin{thm}\label{thm:finite_cone}
 Let  $(T,P,\mu)$ be a BRW with
finitely many cone types and constant mean offspring $m(x)= m
>1.$ We have
\begin{itemize}
    \item[(i)]   $\alpha=0$ if $m\leq 1/\rho(P),$
    \item[(ii)]  $0<\alpha <1$ if $1/\rho(P)< m \leq  1/\widetilde\rho(P),$
    \item[(iii)]  $\alpha=1$ if $m>1/\widetilde\rho(P)$.
\end{itemize}
\end{thm}

\begin{rem}\label{rem:ex_compete}
The example in Theorem \ref{thm:finite_cone} illustrates very well
the two \emph{exponential} effects that compete. The first is the
exponential decay of the return probabilities represented by
$\rho(P_i)$ and the other the exponential growth of the particles
represented by $m.$ If $m$ is smaller that $1/\rho(P_i)$ for all
$i$  the decay of the return probabilities always wins and the
process is transient. In the middle regime where $1/\rho(P_i)< m
<1/\rho(P_j)$ for some $i,j$ the exponential growth may win but if
$m>1/\rho(P_i)$ for all $i$ the exponential growth always wins and
the process is strongly recurrent.
\end{rem}

\subsubsection{Simple random walks on
graphs}\label{subsec:rwgraphs} In order to find conditions for
strong recurrence we inverse the action of connecting graphs at a
common root and split up some given graph in appropriate
subgraphs. In the remaining part of this section we assume for
sake of simplicity the Markov chain to be a simple random walk on
a  graph $G=(V,E),$ where $V=X$ is the vertex set and $E$ is the
set of edges. Keeping in mind that $\rho(P)=\sup_{|F|<\infty}
\rho(P_F),$ compare with equation (\ref{eq:specapprox}), we define
$$\widetilde\rho(P):=\inf_{|\partial F|<\infty} \rho(P_F),$$ where the $\inf$ is over all infinite
irreducible  $F\subset X$  such that the (inner) boundary of $F,$ $\partial
F:=\{x\in F: x\sim F^c\}$ is a finite set. We associate a subset
$F\subset X$ with the induced subgraph $F\subset G$ that has vertex set $F$
and contains all the edges $xy\in E$ with $x,y\in F.$  We can express
$\widetilde\rho(P)$ in terms of transient SRWs on induced subgraphs with
$|\partial F|<\infty.$ For such a graph $F$ we obtain
\begin{eqnarray*}
  p_F^{(n)}(x,y) &=& \P_x(X_n=y, X_i\in F ~\forall i\leq n) \\
   &=& \P_x(X_n=y| X_i\in F~\forall i\leq n)\cdot \P_x(X_i\in F~\forall  i\leq n) \\
   &=& q^{(n)}(x,y)\cdot \P_x(X_i\in F~\forall  i\leq n),
\end{eqnarray*}
where  $q^{(n)}(x,y)$ are the $n$th step probabilities of the SRW
on $F$ with transition kernel  $Q.$ Since $F$ is transient and
$\partial F$ is finite, we have for $x\in F\setminus
\partial F$ that
$$\P_x(X_i\in F~\forall i\leq n)\geq \P_x(X_i\in F~\forall i)>0$$
and hence  $$\limsup_{n\to\infty}
\left(p_F^{(n)}(x,y)\right)^{1/n}=\limsup_{n\to \infty}
\left(q^{(n)}(x,y)\right)^{1/n},\quad \forall x,y\in F.$$
Eventually, we can write
$$\widetilde\rho(P)=\widetilde\rho(G)=\inf_{|\partial F|<\infty} \rho(F),$$ where
the $\inf$ is over all induced connected infinite subgraphs
$F\subset G$ with finite boundaries. In analogy to the proof of
Lemma \ref{lem:rec_common_root} we obtain a necessary condition
for strong recurrence that we conjecture to be sufficient for
graphs with bounded degrees.

\begin{lem}\label{lem:tilderho} The BMC $(G,\mu)$ is not
 strongly recurrent, $\alpha<1,$  if
$$m< \frac1{\widetilde\rho(G)}.$$ If the $\inf$ is attained then
 $m=1/\widetilde\rho(G)$ implies $\alpha<1.$
\end{lem}

\begin{rem}
Lemma \ref{lem:tilderho} holds true for any locally finite graph.
However,  $m>1/\widetilde\rho(P)$ does not imply strong recurrence
in general, see the following Example \ref{ex:str_rec_loc_finite}
and Subsection \ref{sec:outlook} for a more detailed discussion.
\end{rem}

\begin{ex}\label{ex:str_rec_loc_finite}
Consider the following tree $T$  with exploding degrees bearing
copies of $\Z^+$  on each vertex. Let $r$ be the root with degree
$7.$ First, define inductively the \emph{skeleton} of our tree:
$deg(x)=2^{2n+3}-1$ for vertices $x$ with $d(r,x)=n.$ Now, glue on
each vertex a copy of $\Z^+,$ such that in the final tree a vertex
with distance $n$ from the root has degree $2^{2n+3}$ or $2.$ Due
to this construction we have $\rho(T)=\widetilde\rho(T)=1.$
Consider  the BRW $(T,\mu)$ with $\mu_2(x)=1$ for all $x\in T$ and
start the process with one particle in $r.$ The probability that
no copy of $\Z^+$ is visited is at least the probability that the
process lives only on the skeleton and moves always away from the
root:
$$\left(1-\frac18\right)^2\cdot\prod_{n=1}^\infty
\left(1-\frac1{2^{2n+2}}\right)^{2^{n+1}}>0.$$ Hence
the BRW is not strongly recurrent.
\end{ex}

In order to give a sufficient condition for strong recurrence we define
\begin{equation}
\check{\rho}(P):=\limsup_{n\to\infty} \inf_{x\in X}
\rho(P_{B_n(x)})
\end{equation}
and write $\check{\rho}(G)$ for the SRW on $G.$ Here, $B_n(x)$ is
the ball of radius $n$ around $x.$

Notice that this can be seen as a variation of of the spectral
radius since $$\rho(P)=\limsup_{n\to\infty} \sup_{x\in X}
\rho(P_{B_n(x)}).$$

\begin{prop}\label{prop:checkrho_str_rec}
Let $G$ be a graph with bounded geometry. The BRW $(G,\mu)$ with
constant offspring distribution is strongly recurrent if
$$m> \frac1{\check\rho(G)}.$$
\end{prop}
\begin{proof}
 There exists some $n\in \N$ such that for all
$x\in X$ we have $m> 1/ \rho(P_{B_n(x)}).$ We follow the lines of
the proof of Theorem \ref{thm:2} and construct and infinite number
of  supercritical Galton-Watson processes. Observe that since the
maximal degree of $G$ is bounded, there are only a finite number
of different possibilities for the graphs $B_{x,n}.$ Therefore,
the extinction probabilities of the Galton-Watson processes are
bounded away from $1$ and and we can conclude with Lemma
\ref{lem:strongrec}.
\end{proof}

\begin{rem}
The sufficient condition in Proposition
\ref{prop:checkrho_str_rec} is  not necessary for strong
recurrence in general.  Consider the following tree $T$ that is a
combination of $\T_3$ and $\Z:$ Let $r$ be the root with degree
$2.$  The tree $T$ is defined such that $ deg(x)=3$ for all $x$
such that $d(r,x)\in \left[2^{2k}, 2^{2k+1}-1\right]$  and
$deg(x)=2$ for all $x$ such that $d(r,x)\in \left[2^{2k+1},
2^{2k+2}-1\right]$ for  $k\geq 0.$ We have
$\rho(T)=1,~\check\rho(T)=\rho(\T_3)=2 \sqrt{2}/3<1$ and that the
BRW $(T,\mu)$ is strongly recurrent for all $m>1.$ To see the
latter observe that for all $m>1$ there exists some $k$ such that
$p_{\Z}^{(k)}(0,0)\cdot m^k >1,$ where $P_{\Z}$ is the transition
kernel of the SRW on $\Z.$ Thus each part of $\Z$ of length  $k$
constitutes a seed and we conclude with Lemma \ref{lem:strongrec}.
\end{rem}

\subsubsection{Simple random walks on trees}

Let $T$ be a tree of degree bounded by $M\in\N$ and denote $P_T$
for the transition matrix of the SRW on $T$ and $P_{\T_M}$ for the
transition matrix for the SRW on $\T_M,$ the $M$-regular tree. We
consider $T$ to be an infinite subtree of $\T_M.$ One shows by induction on $n:$

\begin{lem}\label{lem:8.1}
$$p_T^{(n)}(x,y)\geq p_{\T_M}^{(n)}(x,y)\quad \forall x,y\in T~\forall
n\in \N.$$
\end{lem}

Since the spectral radius of the SRW on $\T_M$ is
$\rho(\T_M)=\frac{ 2 \sqrt{M-1}}M,$ compare with  Example
\ref{ex:BRWontree}, we immediately obtain a lower bound for the
spectral radius of SRW on trees with bounded degrees.
\begin{lem}\label{lem:spec_tree}
Let $T$ be a tree with degrees bounded by $M.$ Then the simple
random walk on $T$ satisfies
$$ \rho(T)\geq \frac{ 2 \sqrt{M-1}}M.$$
\end{lem}
We obtain the following Corollary of Lemma \ref{lem:spec_tree} and
Proposition \ref{prop:checkrho_str_rec}.
\begin{cor}\label{cor:tree_strongrec}
 Let $T$ be a tree  with maximal degree $M.$
The BRW $(T,\mu)$ with constant offspring distribution  is
strongly recurrent if $m>\frac{ 2 \sqrt{M-1}}M.$
\end{cor}

\begin{rem}
Observe that Lemma \ref{lem:spec_tree} does hold  for general
graphs with degrees bounded by $M,$ compare with  Theorem 11.1 in
\cite{woess}. Therefore, Corollary \ref{cor:tree_strongrec} does
hold true for graphs with degrees bounded by $M$ as well.
\end{rem}

We conclude this chapter with an interesting and illustrative
example, gathered from \cite{woess} (Chapter 9),  where we can
give a complete classification in transience, recurrence and
strong recurrence.\
\\

\noindent
\begin{minipage}{8cm}
\begin{ex}\label{ex:line_tree}
We construct a graph that looks like a rooted $M$-ary tree with a
hair of length $2$ at the root. Let $G_1$ be the tree where each
vertex has degree $M>1,$ with the exception of the root $o,$ which
has degree $M-1.$ As $G_2$ we choose the finite path $[0,1,2].$
The graph $G$ is obtained by identifying $0$  with $o,$ compare
with Subsection \ref{subsec:glue}. The SRW on $G$ is obtained by
setting $\alpha_1=\frac{M-1}M$ and $\alpha_2=\frac1M,$ compare
with Equation (\ref{eq:glueMC}). Let us first consider the case
where $M\geq 5.$ One calculates the spectral radius of the SRW on
$G:$
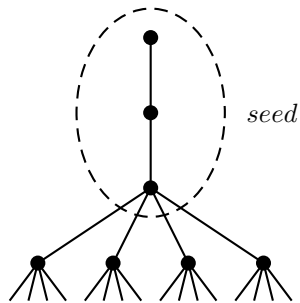
$$\rho(G)=\sqrt{\frac{M-1}{2(M-2)}}.$$
\end{ex}
\end{minipage} \hfill
\begin{minipage}[l]{4cm}
\begin{pspicture}(-2,-2.8)(2,3)
%X_2
\psline{-}(0,0)(0,2) \pscircle*(0,0){0.1}
\pscircle*(0,1){0.1}\pscircle*(0,2){0.1}

%X_2
\pscircle*(0.5,-1){0.1}\pscircle*(1.5,-1){0.1}
\pscircle*(-0.5,-1){0.1}\pscircle*(-1.5,-1){0.1}

\psline{-}(0,0)(0.5,-1)\psline{-}(0,0)(1.5,-1)
\psline{-}(0,0)(-0.5,-1)\psline{-}(0,0)(-1.5,-1)
%\uput[u](0.5,0.5){$\frac12$} \uput[u](-.5,0.5){$\frac12$}

\psline{-}(0.5,-1)(0.875,-1.5)\psline{-}(0.5,-1)(0.625,-1.5)
\psline{-}(0.5,-1)(0.375,-1.5)\psline{-}(0.5,-1)(0.125,-1.5)

\psline{-}(1.5,-1)(1.875,-1.5)\psline{-}(1.5,-1)(1.625,-1.5)
\psline{-}(1.5,-1)(1.375,-1.5)\psline{-}(1.5,-1)(1.125,-1.5)

\psline{-}(-0.5,-1)(-0.875,-1.5)\psline{-}(-0.5,-1)(-0.625,-1.5)
\psline{-}(-0.5,-1)(-0.375,-1.5)\psline{-}(-0.5,-1)(-0.125,-1.5)

\psline{-}(-1.5,-1)(-1.875,-1.5)\psline{-}(-1.5,-1)(-1.625,-1.5)
\psline{-}(-1.5,-1)(-1.375,-1.5)\psline{-}(-1.5,-1)(-1.125,-1.5)

\psellipse[linestyle=dashed](0,1)(1,1.4)\uput[r](1.1,1){$seed$}
\end{pspicture}
\end{minipage}

\noindent Recall that $\vr(G)=\min\{\rho(G_1), \rho(G_2)\}.$ Due
to Lemma \ref{lem:spec_tree} we have $\vr(G)\geq \frac{ 2
\sqrt{M-1}}M.$ Since $\rho(G_1)=\frac{ 2 \sqrt{M-1}}M$ we have
$\vr(G)=\frac{ 2 \sqrt{M-1}}M.$ Notice that
$\vr(G)=\widetilde\rho(G)=\check\rho(G).$ Now, Theorem
\ref{thm:1a}, Lemma \ref{lem:rec_common_root} and the proof of
Theorem \ref{thm:rec_quasi}  yields
\begin{enumerate}
    \item[(i)] $m\leq 1/\rho(G) ~\Longrightarrow$ $(G,\mu)$ is
    transient,
    \item[(ii)] $1/\rho(G)<m\leq 1/\vr(G)~\Longrightarrow$
    $(G,\mu)$ is recurrent,
    \item[(iii)] $m>1/\vr(G)~\Longrightarrow$ $(G,\mu)$ is strongly recurrent.
\end{enumerate}

Observe that in this example the graph $G_2$ can be seen as a seed
that makes the BRW recurrent. The first critical value $1/\rho(G)$
is such that the
 $G_2$ becomes a seed, where the second critical value
$1/\vr(G)$ is such that the branching compensates the drift
induced by the graph $G_1.$ Furthermore, notice that for $M=3,4$
the spectral radius of $G$ is $\rho(G)=\frac{2\sqrt{M-1}}M$ and
recurrence and strong recurrence coincide. Thus in this case, the
branching which is necessary to produce a seed in $G_2$ must be at
least as high as the branching that is needed to compensate the
drift of the SRW on $G_1.$

\subsection{Outlook}\label{sec:outlook} We know that if the
offspring distributions depend on the state, any criterion for
strong recurrence must incorporate more information on the
offspring distributions than the mean. If the offspring
distributions do not depend on the state, we conjecture, compare
with the results obtained  Section \ref{sec:strong}, that there is
a second threshold:

\begin{conj}
Let $(X,P,\mu)$ be a BMC with constant offspring distribution.
Then there exists some $\widetilde m$ such that the BMC is
strongly recurrent if  $m>\widetilde m$ and not strongly recurrent
if $m<\widetilde m.$
\end{conj}

Let us  state the conjecture made in Subsection
\ref{subsec:rwgraphs}. Recall
$$\widetilde\rho(G)=\inf_{|\partial F|<\infty} \rho(F),$$ where the $\inf$
is over all induced connected infinite subgraphs $F\subset G$ with
finite boundaries.
\begin{conj}\label{conj:strongrec1}
Let $G$ be a  graph  with bounded degrees. The BRW $(G,\mu)$ with
constant offspring distribution is
 strongly recurrent if
$$m> \frac1{\widetilde\rho(G)}.$$
\end{conj}
For SRWs on locally finite graphs this is not true, compare with
Example \ref{ex:str_rec_loc_finite}. This example suggests to
consider transient subsets.  Let
$$\widetilde\rho(G,m):=\inf_F \rho(F),$$ where the
$\inf$ is over all irreducible $F\subset G$ where $\partial F$ is
transient  with respect to the BRW $(F,\mu)$. Observe that
$\widetilde\rho(G,m)$ does depend on $m$ since transience is
w.r.t. the BRW. In analogy to the proof of Lemma
\ref{lem:tilderho} we can prove that the BRW $(G,\mu)$ is not
strongly recurrent if $m< 1/\widetilde\rho(G,m).$ We conjecture
that for BRWs   $1/\widetilde\rho(G,m)$ is decisive for strong
recurrence, compare with Lemma \ref{lem:tilderho}.
\begin{conj}\label{conj:strongrec2}
Let $G$ be a  graph. The BRW $(G,\mu)$ with constant offspring
distribution is
 strongly recurrent if
$$m> \frac1{\widetilde\rho(G,m)}.$$
\end{conj}

\section{Positive recurrence}\label{sec:posrec}
An irreducible Markov chain is  called  positive recurrent if the
expected time to return is finite for all possible starting
positions. We generalize this definition to  BMC and say the
process returns to its starting position if the starting position
is hit by at least one particle. If the expected time to return is
finite, we call the BMC  positive recurrent.
\begin{defn}\label{def:posrec_BMC}
A recurrent BMC  is  positive recurrent if
\begin{equation}\label{eq:posrec}
 \E_x T_x <\infty\quad \forall x\in X,
\end{equation}
with $T_x:=\inf\{n>0: \exists i\in\{1,2,\ldots,\eta(n)\}:
x_i(n)=x\}.$ Otherwise it is called  null recurrent.
\end{defn}

In contrast to the question of transience and recurrence of BMC,
it is now also interesting to consider underlying null recurrent
Markov chains and ask whether the corresponding BMC is null or
positive recurrent. For a Markov chain we have that either $\E_x
T_x<\infty$  for all or for none $x\in X.$ This does no longer
hold for BMCs, as we can see in the following Example
\ref{ex:pos_rec}. Furthermore, positive recurrence does depend on more information of the offspring distribution then just the mean offspring.\ \\

\noindent
\begin{minipage}{7.5cm}
\begin{ex}\label{ex:pos_rec}
We consider a random walk on a directed graph with denumerable
many directed cycles of \emph{exploding} length emanating from the
origin $\orig.$ Let $C_i=(c^{(i)}_{0},\cdots,c^{(i)}_{2^i})$ be
cycles of length $2^i$ with $c^{(i)}_{0}=c^{(i)}_{2^i}=\orig$  for
$i\geq 1.$ The origin $\orig$ is the only common vertex of these
cycles, i.e., $c^{(i)}_{k}\neq c^{(j)}_{l},$ for
 $1\leq k< 2^i$ and $ 1\leq l < 2^j~\forall i,j\in\N.$ The transition probabilities $P$ on $X:=\bigcup_i C_i$
are defined as
\end{ex}
\end{minipage}
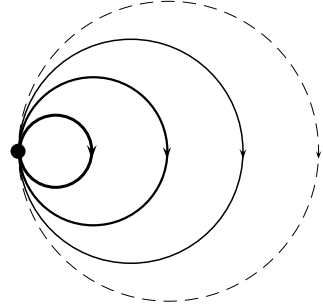
\hfill
  \begin{minipage}[h]{4.1cm}
\psset{linewidth=0.05}
\begin{pspicture}(-0.1,-1.0)(4,1.0)
\pscircle*(-0,0){0.1} \pscircle[linewidth=0.04](0.5,0){0.5}
\pscircle[linewidth=0.03](1,0){1}
\pscircle[linewidth=0.02](1.5,0){1.5}
\pscircle[linewidth=0.01,linestyle=dashed](2,0){2}
\psline[linewidth=0.04]{->}(0.98,0)(0.98,-0.1)
\psline[linewidth=0.03]{->}(1.98,0)(1.98,-0.1)
\psline[linewidth=0.02]{->}(2.99,0)(2.99,-0.1)
\psline[linewidth=0.01]{->}(3.995,0)(3.995,-0.1)
\end{pspicture}\par
\end{minipage}

\noindent
\begin{eqnarray*}
  p(\orig,c^{(i)}_{1}) &:=&  \left(\frac12\right)^i, \quad i\geq 1, \\
  %p(\orig,c^{(i)}_{k}) &:=& 0\quad \forall 1<k\leq 2^i,~i\geq 1 \\
  p(c^{(i)}_{k},c^{(i)}_{k+1}) &:=& 1 \quad \forall 1\leq k<2^i,~i\geq 1.
\end{eqnarray*}

\noindent The Markov chain $(X,P)$  is null recurrent. We consider
the BMC $(X,P,\mu)$  with $\mu_1(\orig)=\mu_3(\orig)=\frac12$ and
$\mu_2(x)=1~ \forall x\neq \orig.$ It is now straightforward to
show that $\E_{\orig}T_{\orig}=\infty$ but
$\E_{c^{(1)}_{1}}T_{c^{(i)}_{1}}<\infty.$ Observe that  the BMC
with the same constant mean offspring $m=2$ but $\mu_2(x)=1$ for
all $x\in X$ is positive recurrent, i.e., $\E_x T_x<\infty$ for
all $x\in X.$

Despite Example \ref{ex:pos_rec} we have under some natural
assumptions on the branching that $\E_x T_x<\infty$ holds either
for all or none $x\in X.$

\begin{lem}\label{lem:pos_rec1}
Let $(X,P,\mu)$ be a BMC and assume that $0<\mu_1(x)<1$ for all $x\in X.$ If $\E_{\orig}
T_{\orig}<\infty$ for some $\orig\in X$ we have
$$ \E_x T_y<\infty \quad \forall x,y\in X.$$
\end{lem}

\begin{proof} We start the BMC in $\orig.$ Let $y\in X$ and choose $k$ such
that $p^{(k)}(\orig,y)>0.$ Since $\mu_1(x)>0~ \forall x$ we have
that with positive probability the total number of particles at time $k$ is $1$ and that this particle is in $y$,  i.e., $\eta(k)=\eta(k,y)=1.$
Hence, $\E_y T_{\orig}<\infty.$ In order to show
$\E_{\orig}T_y<\infty$ we use that $\mu_1(x)<1.$ Let $\tau_i$ be
independent random variables distributed like $T_{\orig}$ under
$\P_{\orig}[\cdot|\eta(1)=1].$  We proceed with a \emph{geometric
waiting time argument}:  We start the process with one particle in
$\orig$ and wait a random time $\tau_i$ until a first particle
returns. This particle splits up in at least two particles with
positive probability $1-\mu_1(\orig).$ One of these particles
starts a new process that returns to $\orig$ after $\tau_2$ time
steps. The remaining particles, if there exists any,  hit $y$
after $k$ time steps with positive probability at least $p^{(k)}(x,y).$
This is repeated until $y$ is hit. Therefore, we obtain with
$q:=(1-\mu_1(\orig))p^{(k)}(\orig,y):$
$$
  \E_{\orig} T_y \leq k+ \sum_{i=1}^\infty \left( (1-q)^{i-1} q\cdot
  \sum_{j=1}^i E \tau_j\right)<\infty,
$$
since $E\tau_i=\E_{\orig}[T_{\orig}|\eta(1)=1]<\infty.$
 Hence, $\E_{\orig} T_y<\infty$ and $\E_y T_{\orig}<\infty$ for all
$y\in X.$ Since $\E_x T_y \leq \E_x T_{\orig} + \E_{\orig}
T_y~\forall x,y\in X$ we are done.
\end{proof}

There is a \emph{branching} analog  to the 2nd criterion of
Foster, compare with Theorem 2.2.3. in \cite{fayolle95}, for
positive recurrence of Markov chains.

\begin{thm}\label{thm:bmc:2Foster}
Let $\orig\in X.$ If there exists a nonnegative function $f$ with
$f(\orig)>0$  such that
\begin{equation}\label{eq:bmc:2Foster}
 Pf(x)\leq \frac{f(x)-\eps}{m(x)} \quad\forall
x\neq \orig~\mbox{for some }\eps>0,
\end{equation}
then $\E_{x}T_{\orig}<\infty$ for all $x\neq \orig.$
\end{thm}

\begin{proof}Let $x\neq \orig.$ We define $$
Q(n):=\sum_{i=1}^{\eta(n)} f(x_i(n))$$ and
\begin{equation}\label{eq:Z(n)}
Z(n):=Q(n\wedge T_{\orig}) + \eps\cdot(n\wedge T_{\orig}).
\end{equation}
We write $\omega(n):=\{x_1(n),\ldots,x_{\eta(n)}(n)\}$ for the
positions of particles at time $n$ and obtain using Equation (\ref{eq:bmc:2Foster}):
$$\E_x[Q(n+1) | \omega(n)=\omega]\leq Q(n) - \eps \eta(n)$$
under $\{T_{\orig}>n\}.$ Hence, under $\{T_{\orig}> n\}$ we have
\begin{eqnarray*}
   \E_x[Z(n+1)|\omega(n)=\omega] &=& \E_x[Q(n+1)+\eps(n+1)]\\
   &\leq&  Q(n)-\eps \eta(n)+ \eps(n+1)\\
    &\leq& Q(n)+\eps n.
 \end{eqnarray*}
 Therefore,  $Z(n)$ is a nonnegative supermartingale. We obtain with
 Equation (\ref{eq:Z(n)})
 $$
\E_x[n\wedge T_{\orig}]\leq \frac{\E_x[Z(n)]}{\eps}\leq
\frac{\E_x[Z(0)]}{\eps} =\frac{ f(x)}{\eps}.$$ Letting
$n\rightarrow\infty$ yields
$$\E_x [T_{\orig}]\leq \frac{f(x)}{\eps}<\infty\quad \forall x\neq \orig.$$
\end{proof}

In general it is not possible to give criteria for the positive
recurrence in terms of the mean offspring $m(x),$ compare with Example
\ref{ex:pos_rec}.  Despite this fact, it turns out that for
\emph{homogeneous} BMC, e.g. quasi-transitive BMC, strong
recurrence and positive recurrence coincide. We refer to
\cite{comets:05} where the asymptotic of the tail of the
distributions of the hitting times are studied even for branching
random walks in random environment.

In the following subsection,we present another method to show positive recurrence of BRW on
$\Z$ using large deviation estimates and the rate of escape of the
BRW.

\subsection{BRW on $\Z$}\label{subsec:pos_rec_Z}
Let us consider an irreducible, transient random walk,
$S_n=\sum_{i=1}^n X_i,$ on $\Z$ with i.i.d. increments $X_i.$ Furthermore, we assume bounded jumps, i.e., $|X_i|\leq d$ for some $d\in\N.$ 
This assumption will be crucial in the proof of Lemmata \ref{lem:ExTx} and \ref{lem:ExTx2} but can be replaced for Theorem \ref{thm:speedBMC1} and Corollary
\ref{cor:speedBMC2} by the assumption that $S_n/n$ satisfies a
large deviation principle. Without loss of generality we assume
that the random walk has drift to the right, i.e., $E[X_i]>0.$ Let
$I(\cdot)$ be the strictly monotone rate function defined by
$$ -I(a)=\lim_{n\to\infty} \frac1n \log \P(S_n\leq a n )
\mbox{ for } a\leq E X_i.$$ We consider the BMC  $(X,P,\mu)$ with
constant offspring distribution $\mu$ with mean offspring $m$ and
denote $M_n$ the leftmost particle at time $n.$

We have the well-known result for the minimal position of a BRW,
compare with \cite{hammersley:74}, \cite{kingman:75},
\cite{biggins:76}, \cite{peres99}, and \cite{menshikov:07} for multi-type BRWs.

\begin{thm}\label{thm:speedBMC1}
$$\liminf_{n\to\infty} \frac{M_n}n=\sup\{s: I(s)\geq \log m\}\quad\P\mbox{-a.s.}$$
\end{thm}

Since $\rho(P)=e^{-I(0)},$ compare with Lemma \ref{lem:spectral_ldp}, we
immediately obtain the following result on the speed of the
leftmost particle.

\begin{cor}\label{cor:speedBMC2} Let $(X,P)$ be a random walk with
bounded jumps on $X=\Z$ and drift to the right. For a BRW
$(X,P,\mu)$ with constant offspring distribution the following
holds true:
\begin{enumerate}
    \item[(i)] If $m>1/\rho(P),$  then
$$\liminf_{n\to\infty} \frac{M_n}n<0\quad \P\mbox{-a.s.}$$
 \item[(ii)] If $m=1/\rho(P),$  then
$$\liminf_{n\to\infty} \frac{M_n}n=0\quad \P\mbox{-a.s.}$$
    \item [(iii)] If $m<1/\rho(P),$ then
$$\liminf_{n\to\infty} \frac{M_n}n>0\quad \P\mbox{-a.s.}.$$
\end{enumerate}
\end{cor}
\begin{rem}
In particular, Corollary \ref{cor:speedBMC2} implies transience if
$m<1/\rho(P)$ and strong recurrence if $m>1/\rho(P).$
\end{rem}

Eventually, we obtain that under the above conditions strong
recurrence implies positive recurrence:
\begin{thm}\label{thm:posrecZ}
Let $(X,P)$ be a random walk with bounded jumps on $X=\Z$ and
drift to the right. The BRW $(X,P,\mu)$ with constant offspring
distribution is positive recurrent if  $m>1/\rho(P).$
\end{thm}

The proof follows from the following Lemmata \ref{lem:ExTx} and
\ref{lem:ExTx2}. Beside $T_x$ we consider the following stopping
time

$$ \widetilde{T}_x:=\inf_{n>0}\{\exists i\in\{1,2,\ldots,\eta(n)\}: x_i(n)\in [x-d,x]\}.$$

\begin{lem}\label{lem:ExTx}
Let $(X,P)$ be a random walk with bounded jumps on $X=\Z$ and
drift to the right. For a BRW $(X,P,\mu)$ with constant offspring
distribution and  $m>1/\rho(P)$  we have
$$ \E_x \widetilde{T}_x<\infty\quad \forall x\in X.$$
\end{lem}
\begin{proof}
We show  the claim for $x=0$ and  write $\widetilde{T}$ for
$\widetilde{T}_0.$ Since
$$\E \widetilde{T} = \sum_{n\geq 0} \P(\widetilde{T}> n)$$ it suffices to study the behavior
of $\P(\widetilde{T}> n)$ for large $n$ and to show that it is
summable. To this end we split the sum into two terms:
\begin{eqnarray*}
 \P(\widetilde{T}>n)&=& \sum_{k=1}^{\lambda n} \P\left(\widetilde{T}>n| \eta(\gamma n)=k\right)
\P\left( \eta(\gamma n)=k\right) \\
  &+& \sum_{k=\lambda
n+1}^\infty \P\left(\widetilde{T}>n| \eta(\gamma n)=k\right)
\P\left( \eta(\gamma n)=k\right)
\end{eqnarray*}
with  $\lambda>0$ and $\gamma\in(0,1)$ to be chosen later.
Here we assume that $\lambda n$ and $\gamma n$ take values in $\N.$ We obtain
\begin{equation}\label{eq:tildeT}
\P(\widetilde{T}>n) \leq \P\left( \eta(\gamma n)\leq \lambda
n\right) + \sum_{k=\lambda n+1}^\infty \P\left(\widetilde{T}>n|
\eta(\gamma n)=k\right).
\end{equation}
In order to estimate the second summand we observe  that at time $\gamma n$ all particles are \emph{at worst} at position $\gamma n d$ and obtain
\begin{eqnarray*}
  \P(\widetilde{T}>t|\eta(\gamma n)=k) &\leq& \P\left(M_{(1-\gamma)n} >
  -\gamma n d\right)^k    \\
   &=& \P\left( \frac{ M_{(1-\gamma)n}}{(1-\gamma)n} >
   -\frac{\gamma d}{1-\gamma}\right)^k.
\end{eqnarray*}
Due to Corollary \ref{cor:speedBMC2} we have $$ \liminf_n \frac{
M_{(1-\gamma)n}}{(1-\gamma) n}=s$$ for some $s<0.$ Now, we choose
$\gamma$  such that
$$ -\frac{\gamma d}{1-\gamma}>s.$$
Hence, there exists  $\theta<1$ with
$$\P\left( \frac{ M_{(1-\gamma)n}}{(1-\gamma)n} > -\frac{\gamma
d}{1-\gamma}\right)\leq \theta<1,$$ for sufficiently large $n.$
Therefore,

\begin{equation}\label{eq:lem:ExTX}
 \P(\widetilde{T}>n|\eta(\gamma n)=k)\leq \theta^k
\end{equation} and the second summand in Equation
(\ref{eq:tildeT}) is bounded by
$\theta^{\lambda(n+1)}/(1-\theta).$

  It remains to bound the first term in Equation
(\ref{eq:tildeT}). To do this we do not consider the whole BRW but
focus on the induced random walk. Denote $Y_n$ the number of times
the labelled particle is not the only offspring of its ancestor.
In other words, when we think about the process where particles
live forever and produce offspring according to
$\tilde\mu_{i-1}=\mu_i~i\geq 1,$ then $Y_n$ is just the number of
offspring  of the starting particle at time $n.$ Hence,  $Y_n\sim
Bin(n,p),$ where $p:=\sum_{i=2}^\infty\mu_i>0.$ Observe that  a
Large Deviation Principle holds for $Y_n$, i.e., $P(Y_n\leq a n)$
decays exponentially fast for $a<p.$ Due to the definition of
$Y_n$ we have  $\eta(\gamma n)\geq Y_{\gamma n}$ and obtain with
$l:=\gamma n$ that
$$ P(\eta(\gamma n)\leq \lambda n)\leq P(Y_{\gamma n}\leq \lambda
n)= P(Y_l\leq \lambda \gamma^{-1} l).$$ The last term decays
exponentially fast for $\lambda\gamma^{-1}<p.$ Therefore, choosing
$\lambda<p\gamma,$ we obtain   that $\P\left( \eta(\gamma n)\leq
\lambda n\right)$ decays exponentially fast.
\end{proof}

\begin{lem}\label{lem:ExTx2} Under the assumptions of Lemma \ref{lem:ExTx}
we have
$$ \E_x \widetilde{T}_x<\infty ~\Longrightarrow \E_x T_x<\infty. $$
\end{lem}

\begin{proof}
For simplicity let $x=0$ and write $\widetilde{T}$ and $T$ for
$\widetilde{T}_0$ and $T_0,$ respectively. Analogous to the proof
of Lemma \ref{lem:ExTx}, it is proven that $\E_y
\widetilde{T}<\infty~\forall y\in X.$ We start the BRW in $0$ and
consider the random time $\tau_1$ when $[-d,0]$ is hit for the
first time by some particle. We pick one particle of those being
in $[-d,0]$ at time $\tau_1$ and consider the BRW originating from
this particle.  Due to the irreducibility, there exists $k\in\N$
and $q>0$ such that $\forall y\in [-d,0]$ we have
$p^{(l)}(y,0)\geq q>0$ for some $l\leq k.$ Hence, $0$ is visited
by some particle up to time $k$ with probability at least $q.$ If
$0$ is not hit after $k$ time steps  we consider a BRW starting in
some occupied position in $[-d(k+1),d(k+1)]$ and wait the random
time $\tau_2$ until  $[-d,0]$ is hit by some particle. This is
repeated until $0$ is hit at the random time $W.$ Since $\E_x
\widetilde{T}<\infty$ for all $x\in X,$ there exists some $C>0$
such that $\E_x \widetilde{T}\leq C ~\forall x\in [-d(k+1),
d(k+1)].$ We conclude with
\begin{eqnarray*}
  \E T\leq \E W &\leq& \E\left[ (\tau_1+k)q+ (1-q)q (\tau_1+\tau_2+k)+\cdots \right]\\
   &\leq& \sum_{i=1}^\infty (iC+k) (1-q)^{i-1} q < \infty.
\end{eqnarray*}
\end{proof}

\def\cprime{$'$}

\end{document}